\documentclass[journal]{IEEEtran}
\usepackage[cmex10]{amsmath}
\usepackage[pdftex]{graphicx}
\usepackage{pdfsync}
\usepackage{amssymb}
\usepackage{color}
\usepackage{amsthm}
\usepackage{graphicx}
\usepackage{empheq}
\usepackage{hyperref}
\usepackage{array}
\usepackage{cite}

\newcommand{\cqfd}{%
\mbox{}%
\nolinebreak%
\hfill%
\rule{2mm}{2mm}%
\medbreak%
\par%
}
\newtheorem{assum}{Assumption}

\newtheorem{remark}{Remark}
\newtheorem{lem}{Lemma}
\newcommand{\sgn}{\mathop{\mathrm{sgn}}}
\def\f{\frac}

\def\i1n{i=1,\cdots,n}
\def\j1n{j=1,\cdots,n}
\def\ij1n{i,j=1,\cdots,n}

\newcommand{\be}{\begin{equation}}
\newcommand{\ee}{\end{equation}}
\newcommand{\beq}{\begin{equation*}}
\newcommand{\eeq}{\end{equation*}}

\newtheorem{thm}{Theorem}

\begin{document}

\title{ Time- and State-Dependent Input Delay-\\Compensated Bang-Bang Control of a \\Screw Extruder for 3D Printing
}

\author{
Mamadou~Diagne$^{1},$\thanks{$^{1}$Department of Mechanical and Aerospace Engineering,
University of California, San Diego,
La Jolla, CA, 92093, USA.
 E-mail: {\tt mdiagne@ucsd.edu \ \textrm{and} \ krstic@ucsd.edu}.}
\ Nikolaos Bekiaris-Liberis$^{2},$\thanks{$^{2}$Department of Production Engineering and Management, Technical University of Crete, Chania, 73100, Greece.
E-mail: {\tt nikos.bekiaris@gmail.com }.} and 
\
\ Miroslav~Krstic$^{1}$\thanks{ 
 {\tt }.}
}


\maketitle
\begin{abstract}
In this paper a delay-compensated Bang-Bang control design methodology for the control of the nozzle output flow rate of screw-extruder-based 3D printing processes is developed. The presented application has a great potential to move beyond the most commonly used processes such as Fused Deposition Modeling (FDM) and Syringe Based Extrusion (SBE), improving the build speed and the 3D parts accuracy. A geometrical decomposition of the screw extruder in a partially and a fully filled regions (PFZ and FFZ) allows to describe the material convection in the extruder chamber by a 1D hyperbolic Partial Differential Equation (PDE) coupled with an Ordinary Differential Equation (ODE). After solving the hyperbolic PDE by the Method of Characteristics (MC), the coupled PDE-ODE's system is transformed into a nonlinear state-dependent input delay system. The aforementioned delay system is extended to the non-isothermal case with the consideration of periodic fluctuations acting on the material's convection speed, which represent the effect of viscosity variations due to temperature changes in the extruder chamber, resulting to a nonlinear system with an input delay that simultaneously depends on the state and the time variable. Global Exponential Stability (GES) of the nonlinear delay-free plant is established under a piecewise exponential feedback controller that is designed. By combining the nominal, piecewise exponential feedback controller with nonlinear predictor feedback the compensation of the time- and state-dependent input delay of the extruder model is achieved. Global Asymptotic Stability (GAS) of the closed-loop system under the Bang-Bang predictor feedback control law is established when certain conditions, which are easy to verify, related to the extruder design and the material properties, as well as to the magnitude and frequency of the material’s transport speed variations, are satisfied. Several simulations results are presented to illustrate the effectiveness of the proposed control design.

\end{abstract}

 \section{Introduction}
\begin{figure}[!h]
\setlength{\unitlength}{1.5pt}
\begin{minipage}[t]{0.5\textwidth}
\begin{center}
\includegraphics[width=100\unitlength]{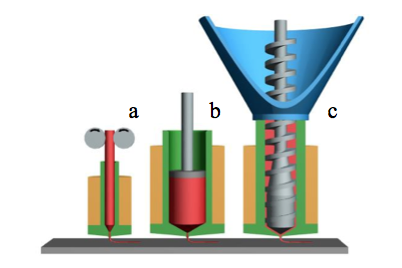}

\caption{\label{fig3}(a) FDM,  (b) SBE,  (c) SE processes  \cite{valkenaers2013}.}
\end{center}
\end{minipage}
\end{figure}

 Additive Manufacturing (AM)  has a  promising future and demonstrates its effectiveness 
in  various  applications involving tissue engineering  \cite{Mironov2003, Billiet2014}, chemical engineering \cite{Dragone2013},  thermoplastics \cite{valkenaers2013},  metal \cite{ADMA2013}  and  ceramic \cite{JBM2005} material's fabrication. Functional 3D objects with complex geometrical shape can be produced in a short time without the needs of tools thanks to the Computer Aid Design (CAD) that drastically reduces the products development procedure. Currently,  the most popular  plastics 3D printers are based on FDM \cite{Widmer1998, Zein2002, JBM2001} and  SBE \cite{wei2013, valkenaers2013} technologies (Fig.\ref{fig3}).

 In these processes, biodegradable polymers are transported, heated and pressurized in an extruder chamber  before being dropped  on a platform,  one horizontal thin layer at a time, until the complete 3D part is built such that it closely resembles the original CAD model. One of the crucial point that is not commonly addressed in the existing literature of extrusion-based 3D printing is controlling of the start and stop of extrusion-on-demand.   A hybrid extrusion force-velocity modeling and tracking control  for the fabrication of functionally graded material parts is developed in \cite{Brad2013} using a first order differential equation that  describes  the plunger dynamic in  a SBE process. Some extents of that approach are proposed by  \cite{Li2011} with a robust tracking of the extrusion force  to recover constant  flow disturbances whereas  \cite{Xiyue2008} considers  an unknown transfer function gain with an adaptive control strategy.  Several issues regarding on  FDM are discussed in  \cite{Michael2007, Satish2007}  and references therein, including the  potential clogging due to agglomerate formation at the nozzle, appearance of bubbles, density inhomogeneity,  tracking of short time-scale process variations,  and prediction of  anomalies such as material overflow and underflow for diverse applications. Thermal control  is left out of most prior studies which are essentially  based on empirical models.

In this paper, we are interested in the flow  control issues related to the  recent advances  of   3D printing technology  for which a Screw Extrusion (SE) process  is utilized.  In SE, the rotating screw allows a continuous feeding mechanism and generates a sufficiently high pressure in the extruder chamber, increasing, as a result, the printing speed. In addition, the screw motion extends the mixing capabilities of the system, and thereby, reduces drastically the risk of potential clogging at the nozzle while improving the homogeneity of the extruded filament \cite{valkenaers2013,Liu06}. The SE process with granular material moves beyond the restrictions of FDM and does not require filament-shaped raw materials to operate. Consequently, it enables the processing of a broader range of raw materials and permits an easy recycling of wasted plastic during extrusion \cite{valkenaers2013,Liu06}. In such processes, the need to control the start and stop of the extrusion process on demand calls for advanced control methodologies that are capable of enhancing the final product's quality in an industrial level.  Even if  experimental results demonstrate  the effectiveness of SE  \cite{Liu06}, \cite{valkenaers2013},  the challenging control problems arising in such applications are actually  poorly investigated.

In the present article, a generic and dynamical model of a homogeneous melt SE process derived from mass and momentum balance laws \cite{ Diagne13, diagne_PhD} is used for the design of a delay-compensated “Bang-Bang” controller which permits a fast and accurate control of the flow at the nozzle output. The model consists of a 1D Partial Differential Equation (PDE) that is defined on a time-varying spatial domain whose dynamics obey to an Ordinary Differential Equation (ODE). The transformation of the coupled PDE-ODE system into a state-dependent input delay system, which describes the dynamics of the material convection in the extruder chamber, is achieved after solving the PDE by the Method of Characteristics (MC)  \cite{Diagne13,Diagne2015}. In order to also account for potential periodic fluctuations of the material’s transport speed when processing granular pellets \cite{Roberts86}, due to the thermal energy that is supplied into the system from the heater of the extruder and due to the mechanical shearing effect by the rotation of the screw,  the state-dependent input delay model is extended  to a nonlinear system with an input delay that depends simultaneously on the state and the time variable (see \cite{Bekiaris13,BEK13} and \cite{bekiaris2012nonlin} for the treatment of systems with time- and state-dependent delays).

In \cite{Diagne2015}, a delay-compensated Bang-Bang control law is developed for the control of the nozzle output flow rate of an isothermal screw extrusion process, achieving GES of the delay-free plant at any given setpoint. By combining the nominal, piecewise exponential feedback controller  \cite{Diagne2015} with nonlinear predictor feedback, which is extended from the state-dependent input delay case \cite{Bekiaris13} to the case in which the vector field and the delay function depend explicitly on time, the compensation of the time- and state-dependent input delay of the non-isothermal screw extrusion model is achieved. GAS of the closed-loop system under the delay-compensated Bang-Bang controller is established when certain conditions, related to the extruder design and the material properties, as well as to the periodic fluctuations, are satisfied. Several simulations results are presented including the case in which there is uncertainty in the value of the periodic variations of the material's transport speed.

This paper is organized as follows: The screw extruder mechanisms and the bi-zone model of the extruder consisting of the transport PDE coupled with the ODE for the moving interface is discussed in Section \ref{model}.  In Section  \ref{transformation}, the transformation of the coupled PDE-ODE system into a state-dependent input delay system by computing the PDE’s solution by the MC is presented and it is then extended to a nonlinear system with a time- and state-dependent input delay. The control of the delay-free plant with a piecewise exponential Bang-Bang-like control law is described in Section \ref{stabilization}. In Section \ref{sectionpredict}, we design the predictor feedback control law for nonlinear systems with time- and state-dependent delay acting on the input. The application of the predictor feedback control law to the screw extruder model is presented in Section \ref{application}. The paper ends with simulations, including a discussion on the robustness properties of a state-dependent input delay compensator to time- varying perturbations acting on the vector field and the delay function, in Section \ref{simulation}.

\section{3D Printing Based on Single-Screw Extruders}\label{model}
\subsection{Extrusion process description and structural decomposition of the extruder into a partially and a fully filled zone}
\begin{figure}[!h]
\setlength{\unitlength}{1.5pt}
\begin{minipage}[t]{0.5\textwidth}
\begin{center}
\includegraphics[height=100\unitlength]{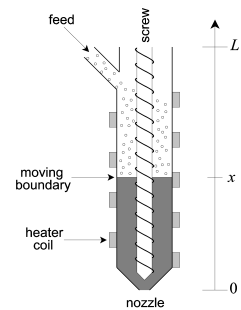}
\caption{\label{Bizone}Bi-zone model of a screw extruder.}
\end{center}
\end{minipage}
\end{figure}

 A screw extruder is divided into one or several conveying zones (transport zones), melting zones (for material fusion) and mixing zones in which the extruded  melt is submitted to high pressure, before its eviction through the nozzle  \cite{Tadmor74, KIM001, KIM002, JAN01, JAN03, CHOU04, KULSH92, Diagne13, diagne_PhD}.  The net flow rate at the extruder nozzle depends  mainly on  the material flow in the longitidunal  direction given by 1D heat and mass transport equations  \cite{Booy80, Booy78}.  Another particularity of these processes  is that they  can be divided in geometric regions which are partially and fully  filled called PFZ and FFZ, respectively (Fig. \ref{Bizone}).  The  PFZ which is submitted to an atmospheric pressure is a conveying region and   the flow in the FFZ  is determined by     the pressure gradient   building-up in that region \textcolor{black}{ due to the nozzle resistance}.  These two zones are coupled by an   interface   which moves according to the volume of material accumulated in the FFZ. Basically, the moving  interface is located at the point where the pressure gradient passes from zero to a non null value. 
\subsection{Mass and momentum balance of an extrusion process}
\subsubsection{Mass  balance of the PFZ }
The PFZ is defined on the time-varying spatial interval $(x(t),L)$, $x(t)$, being the length of the FFZ  and $L$ the extruder length  (Fig. \ref{Bizone}). The mass balance in this area can be expressed using  the fraction of the effective volume  between a  screw element and the barrel ($V_{\rm eff}$)  which is occupied by the extruded material, namely,  the filling ratio  $u$. Considering  an incompressible homogeneous mixture with constant density $\rho_{0}$ and  viscosity $\eta$, the following mass conservation equation is   deduced
\begin{eqnarray}
\partial_{t}u(z,t)&=&\xi N_0\partial_{z}u(z,t),\;(t,z)\in(\mathbb{R}^{+},(x(t), L))\label{transport-fillin-ratio}\\
u(L,t)&=& U(t),\label{boundary}
\end{eqnarray}
where $N_0$ is the constant  screw speed and  $\xi$ the uniform pitch of the screw.   The boundary condition $u(L,t)$ is  defined assuming the continuity of  the  flow at the inlet $\{z=L\}$ 
\begin{eqnarray}\label{flow}
U(t)=\frac{F_{\rm in}(t)}{\rho_0 N_0V_{\rm eff}},
\end{eqnarray}
where  $F_{\rm in}(t)$ is  the feeding rate. Physically, the term $\rho_0 N_0V_{\rm eff}$ in  (\ref{flow}) is  the maximum pumping capacity of the screw.

\subsubsection{Momentum balance of the FFZ}

The FFZ whose filling ratio is equal to one is defined on  the spatial domain
$(0,x(t))$, where the coordinate $\{z=0\}$ is the extruder's  end.
 The FFZ  flow depends on  \textcolor{black}{the pressure gradient that appears in this region, resulting  to  backward or  forward flow}.  The momentum balance which   is derived from Navier-Stokes equations under stationary conditions yields  the pressure gradient 
\begin{eqnarray}
\partial_{z}P(z,t)=-\eta
\frac{\rho_{0}V_{\rm eff}N_0-F_{d}(t)}{B\rho_{0}},\label{pressure-flow}
\end{eqnarray}
for all $(t,z)\in(\mathbb{R}^{+}, (0, x(t)))$, where $B$ is  a coefficient of pressure flow.
The net  flow rate $F_{d}(t)$, in the case of a Poiseuille flow  is expressed with the help of the nozzle conductance $K_{d}$, the viscosity $\eta$,  and the pressure at the nozzle $P(0,t)$ as
\begin{align}
 & 
  \begin{cases}
F_{d}(t) &=\frac{K_{d}}{\eta}\Delta P(t),\\
\Delta P(t) &=P(0,t)-P_{0}.\label{Die-flow}
  \end{cases}
\end{align}

\subsubsection{Mass balance of the FFZ }

The FFZ mass balance leads to an ODE which describes the time evolution of its  length.  This length denoted by $x(t)$  determines  the location of the small transfer region \textcolor{black}{that is  assimilated to the point}  at which  the pressure changes from the  atmospheric  pressure  $P_0$ to a different value \cite{JAN01, JAN03, KULSH92, CHIN01, diagne_PhD} 
\begin{eqnarray}
\frac{dx(t)}{dt}=\frac{\rho_{0}\xi N_0V_{\rm eff}u(x(t),t) -F_{d}(t)  }{\rho_{0}S_{\rm eff}(1-u(x(t),t))},\label{eq-interface-final}
\end{eqnarray}
where $S_{\rm eff}$ is the available section and $V_{\rm eff}=\xi S_{\rm eff}$. 
\subsubsection{Coupling relations at the PFZ-FFZ interface  }
The coupling condition is formulated imposing the pressure continuity at  the spatial coordinate $x(t)$
\begin{eqnarray}
P(x^-,t)=P(x^+,t)=P_{0}.
\end{eqnarray}
Integrating the pressure gradient  equation (\ref{pressure-flow}), the net flow rate defined in \eqref{Die-flow} is written as 
 \begin{equation}
F_{d}(t)=\frac{K_{d}V_{eff}N_0\rho_{0}x(t)}{B\rho_{0}+K_{d}x(t)}.
\label{pressure-conti}
\end{equation}
\textcolor{black}{Substituting \eqref{Die-flow} in \eqref{eq-interface-final} and using  (\ref{pressure-conti}),   equation (\ref{eq-interface-final}) for the length of the FFZ  is written as}
\begin{eqnarray}
\frac{dx(t)}{dt}&=&  -\xi N_0 \frac{K_d  x(t)-\left(B\rho_0 +K_d x(t)\right) u(x(t),t)}{\left(B\rho_0 +K_dx(t)\right) \left(1-u (x(t),t)\right)}.\label{delay-sys-1}
\end{eqnarray}


 \section{From  Mass Balance Equations of the Extruder to a Delay System} \label{transformation}
\subsection{Isothermal  delay system model}
The  bi-zone model (\ref{transport-fillin-ratio}),  (\ref{boundary}), and
(\ref{delay-sys-1}) can be reduced to a nonlinear state dependent-input delay system \cite{Bekiaris13}. The characteristic solutions of    (\ref{transport-fillin-ratio}) with respect to the boundary condition \eqref{boundary}   are
\begin{eqnarray}
u(z,t) = U\left(t-\frac{L-z}{\xi N_0}\right).\label{init}
\end{eqnarray}
Substituting \eqref{init} into   \eqref{delay-sys-1}, we derive the  following nonlinear system
\begin{align}
\dot x(t)=&  \xi N_0 \left[-\frac{ K_d  x(t)}{\left(B \rho_0 +K_dx(t)\right) \left(1-U\left(t- D_s(x(t))\right)\right)} \right. \nonumber \\&+ \left. \frac{ U \left(t- D_s(x(t))\right)}{ \left(1-U\left(t- D_s(x(t))\right)\right)}\right], \quad U(t) \in [0,1)\label{delay-sys}
\end{align}
The state-dependent input delay function  is denoted as
\be
  D_s(x(t))= \frac{L-x(t)}{\xi N_0}.
\ee 
A detailed derivation of the ODE  \eqref{delay-sys}  and an extensive description of the screw-extruder model for 3D printing  is given in \cite{Diagne2015}.
\subsection{Delay system representation in a non-isothermal case}
In this section,  we propose the following extension of the state-dependent delay model (\ref{delay-sys}) to account for the viscosity changes due to the temperature variations occurring in the physical system when processing plastic pellets 
\begin{align}
\dot x(t)=&  c(t) \left[-\frac{ \theta_2  x(t)}{\left(1 +\theta_2x(t)\right) \left(1-U\left(t- D(t,x(t))\right)\right)} \right. \nonumber \\&+ \left.\frac{ U \left(t- D(x(t))\right)}{ \left(1-U\left(t- D(t,x(t))\right)\right)}\right], \quad U(t) \in [0,1),\label{delay-sys-perturb}
\end{align}
where 
\begin{align}
\theta_1&= \xi N_0,
 \label{theta-1-theo}
\\
\theta_2&= \f{K_d}{B \rho_0},
 \label{theta-2-theo}
\end{align}
and,
\begin{align}
  D(t,x(t))&= \frac{L-x(t)}{c(t)}\label{delay-time},\\
c(t)&=\theta_1 (1+\epsilon \cos(\omega t))\label{delay-time-fluc},
\end{align}
where,  $\epsilon<1$ is  a positive constant  and  $\omega$ is the mean value of the  angular frequency of the periodic fluctuations.
Our choice for the non-isothermal model \eqref{delay-sys-perturb} is motivated by the fact that  the   expansion of granular material    into a plastic state due to the  thermal effect leads implicitly to periodic fluctuations of the convection   speed, namely, $\theta_1$  due to  the viscosity variations  \cite{Roberts86}. Moreover,  some  nozzle ``instabilities" phenomena  may  appear as short periodic distortions of the extrudate, due to the viscoelastic properties of the fluid, with magnitude smaller than one \cite{Roberts86}.

\section{Control of the Delay-Free System with a ``Bang-Bang" Control Law}\label{stabilization}
\subsection{Open-loop stability}
 The starting point  of  the  delay system controller design  consists of the construction of a  nonlinear  control law that stabilizes the   delay-free system 
\begin{align}
\dot x(t)=&  c(t) \left[-\frac{ \theta_2  x(t)}{\left(1 +\theta_2x(t)\right) \left(1-U\left(t\right)\right)} \right. \nonumber \\&+ \left. \frac{ U \left(t\right)}{ \left(1-U\left(t\right)\right)}\right], \quad U(t) \in [0,1)\label{nondelay-sys}.
\end{align} 
For $\epsilon<1$, the time-varying \textcolor{black}{speed of the material transport} $c(t)$ is strictly positive and   the open-loop stabilizing control law of  the delay-free plant  (\ref{nondelay-sys})  is given by
 \begin{equation}
v(x^*)=\frac{\theta_2 x^*}{1+\theta_2x^*}, \quad \forall x^*\in[0,L),
\label{equilibrium}
\end{equation}
for the physical parameters of the extruder satisfying $\theta_2 L<1$.
  This  statement is directly derived considering the Lyapunov function  $V= |e(t)|$, where $e(t)=x(t)-x^*$. 
\subsection{``Bang-Bang" controller design  with piecewise exponential functions}

For the feedback stabilization of (\ref{nondelay-sys}),  we consider   two exponential functions  \cite{Diagne2015}:
\begin{itemize}
\item for $ x(t)\leq x^*$, a  left-exponential function  
\be
\hspace{-0.5cm}v_{\rm l}(x,x^*)=v(x^*)+ (v_{\rm max}-v(x^*)) \frac{1-e^{a_{\rm l}(x^*)(x-x^*)}}{1-e^{-a_{\rm l}(x^*)(x^*)}}\label{left-control},
\ee
where  $a_{\rm l}(x^*)>0$ is  the gain of the left exponential control law.
The function \eqref{left-control}  takes  values  in   $[v(x^*),v_{\rm max}]$, where   $v_{\rm max}<1$   is    the maximal value of  the inlet  filling ratio, namely the maximal feeding capacity of the extruder.  Therefore,  $v_{\rm l}(0)=v_{\rm max}$,  allows to set  the  inlet flow at its maximum capacity for a  rapid refill action when the extruder is empty. 

\item  for $ x(t)\geq x^*$, a right-exponential function 
\begin{eqnarray}
v_{\rm r}(x,x^*)=v(x^*)- v(x^*) \frac{1-e^{-a_{\rm r}(x^*)(x-x^*)}}{1-e^{-a_{\rm r}(x^*)(L-x^*)}}\label{right-control},
\end{eqnarray}
where $a_{\rm r}(x^*)>0$ is  the gain of the right exponential  control law.
The function \eqref{right-control}   belongs  into the interval   $[0\,,\,v(x^*)]$ and   the  control action stops radically the flow when the extruder is completely filled, namely,  $v_{\rm r}(L)=0$.  
\end{itemize}
\subsection{Extension of the ``Bang-Bang"  control law on the whole domain $(0,L)$}  
Next, we introduce  the  characteristic function of the domains \textcolor{black}{ $[0,\ x^*]$ and $[x^*, \ L]$}  and write the extended control law as 
\begin{eqnarray}
v(x,x^*) = v_{\rm l}(x,x^*)h(x^*-x) + v_{\rm r}(x,x^*)h(x-x^*)\label{extend-control},
\end{eqnarray}
where $h$ is  the Heaviside function.

 A continous slope function   at  the setpoint  $x^*$ denoted by $S(x^*)$  is imposed  to  extend  the left and the right exponential  controllers (\ref{left-control}) and  (\ref{right-control}), respectively    into  the  differentiable  piecewise exponential feedback law (\ref{extend-control}).  The   slope function is  defined as  $S(x)=- \frac{dv(x,x^*)}{dx}$ (the minus sign is conventional).  More precisely, the key point of the design is to define a free parameter that may be specified  by the user    as the value of slope function at the equilibrium  $S(x^*)$, under some restrictions that will be emphasize in this section. It is clear that, equiting the assigned  value  $S(x^*)$   to both left and right slope functions  of (\ref{left-control}) and (\ref{right-control}), we can easily  derive the following  relations
\begin{eqnarray}
S(x^*)&=&  \frac{a_{\rm l}(x^*)(v_{max}-v(x^*))}{1-e^{-a_{\rm l}(x^*) x^*}},\label{left-trans-eq}\\
S(x^*)&=& \frac{a_{\rm r}(x^*) v(x^*) }{1-e^{-a_{\rm r}(x^*) (L-x^*)}}.\label{right-trans-eq}
\end{eqnarray}
  The   equations (\ref{left-trans-eq}) and (\ref{right-trans-eq}) are both  transcendental and  admit  numerical solutions namely  the suitable  exponential parameters needed to the left and to the right of the setpoint  to achieve the  differentiability of the controller (\ref{extend-control}).  These solutions  $a_{\rm l}(x^*)>0$ (respectively, $a_{\rm r}(x^*)>0$) exist if  the linear and exponential functions of $a_{\rm l}(x^*)$ (respectively, $a_{\rm r}(x^*)$) have a strictly  positive intersection. Consequently, the  desired slope function $S(x^*)$  should  be  above  some minimum value denoted $S_{\rm min}(x^*)$ ,  for any given equilibrium in the physical domain  $(0,L)$.   More precisely,  the two equations in  (\ref{left-trans-eq}) have   strictly positive solutions  if at the origin  ($a_{\rm l}=0$ and  $a_{\rm r}=0$),  the slope of  their linear part   is less than the slope of their exponential part respectively. 
\begin{itemize}
\item For the left exponential slope function (\ref{left-trans-eq}), we define the linear and the exponential functions
\begin{align}
 & 
  \begin{cases}
\psi_{\rm l}(a_{\rm l}(x^*)) &= a_{\rm l}(x^*)(v_{\rm max}-v(x^*))\\
\phi_{\rm l}(a_{\rm l}(x^*)) & =S(x^*)(1-e^{-a_{\rm l}(x^*)x^*}),
  \end{cases} \label{left-trans-eq-decomposition}
\end{align}
a solution of  (\ref{left-trans-eq}) should  satisfy
\begin{align}
 & 
  \begin{cases}
\psi_{\rm l}(a_{\rm l}) & =  \phi_{\rm l}(a_{\rm l}),\\
\frac{d\psi_{\rm l}(0)}{da_{\rm l}} & < \frac{d\phi_{\rm l}(0)}{da_{\rm l}}.
  \end{cases} \label{left-trans-eq-decomposition-deriv}
\end{align}
It follows that for all  $x \in (0,x^*)$, $S(x^*)$  should satisfiy the  inequality 
\begin{eqnarray}
S(x^*) >\frac{v_{max}-v(x^*)}{x^*}. \label{condition-left} 
\end{eqnarray}
\item Decomposing the  right exponential slope  function (\ref{right-trans-eq}) into
\begin{align}
 & 
  \begin{cases}
\psi_{\rm r}(a_{\rm r}(x^*)) &= a_{\rm r}(x^*)v(x^*)\\
 \phi_{\rm r}(a_{\rm r}(x^*))& =S(x^*) (1-e^{-a_{\rm r}(x^*) (L-x^*)}),
  \end{cases}  \label{right-trans-eq-decomposition}
\end{align}
we deduce that a solution of  (\ref{right-trans-eq}) should  satisfy:
\begin{align}
 & 
  \begin{cases}
\psi_{\rm r}(a_{\rm r}(x^*)) &=  \phi_{\rm r}(a_{\rm r}(x^*)),\\
\frac{d\psi_{\rm r}(0)}{da_{\rm r}} & < \frac{d\phi_{\rm r}(0)}{da_{\rm r}}.
  \end{cases}  \label{right-trans-eq-decomposition-deriv}
\end{align}
\end{itemize}
Hence,  for all  $x \in (x^*, L)$
\begin{eqnarray}
S(x^*) >\frac{v(x^*)}{L-x^*} \label{condition-left} 
\end{eqnarray}

\begin{figure}[!ht]
\setlength{\unitlength}{1pt}
\begin{minipage}[t]{1\textwidth}
\includegraphics[height=120\unitlength]{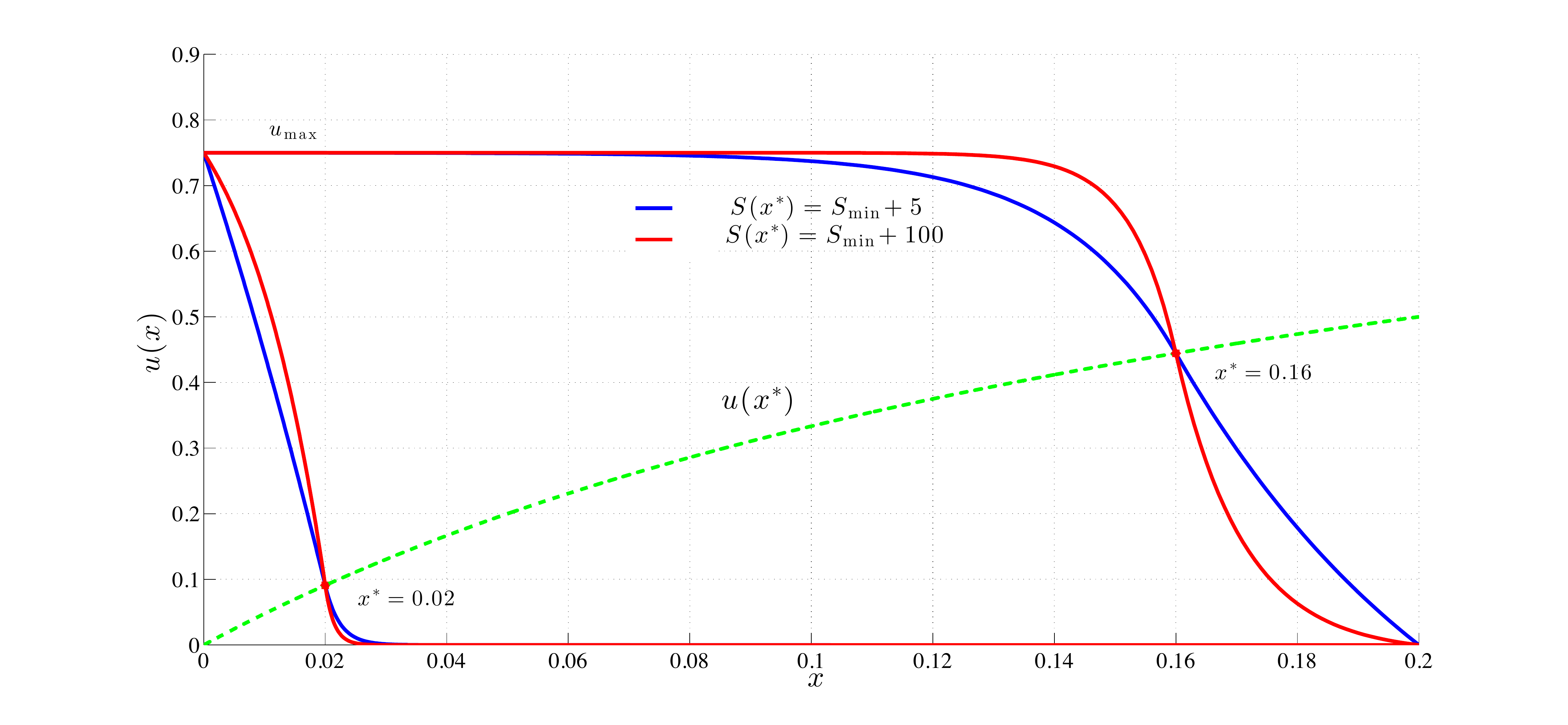}
 \caption{\label{Aggressivness} Control function for different $S_{\rm min}$ and setpoints. }
\end{minipage}
\end{figure}

Finally, the  minimal value of the setpoint  slope  $S_{\rm min}(x^*)$ above which the  gains $a_{\rm l}(x^*)$ and $a_{\rm r}(x^*)$  ensure the differentiability  of the  extended control  law (\ref{extend-control}) on $(0,L)$ is given by
\begin{align}
 S_{\rm min}(x^*)&= \frac{ x^*} {\f{1}{\theta_2x^*}+1} \max \left\{ \frac{v_{\rm max}\left(1+\frac{1}{\theta_2x^*}\right) -1} {x^*}; \frac{1}{L-x^*}\right \}
\label{condition-domain} 
\end{align}
The speed of the controller or its ``agressivness"  increases with the  rise   of the setpoint slope $S(x^*)$. As it is illustrated in   Fig. \ref{Aggressivness}, with the characteristics of the control law for the setpoints $x^*=0.02m$ and $x^*=0.16 m$ with different values of the setpoint slope value $S(x^*)$.

\begin{thm}
For any setpoint $x^* \in [0,L)$ and for any  choosen  setpoint slope  $S(x^*)\in \mathbb{R}$    satisfying $S(x^*)\geq S_{\rm min}(x^*)$, where $S_{\rm min}(x^*)$  is given by (\ref{condition-domain}), taking the control gains $(a_ {\rm l}(x^*),a_ {\rm r}(x^*)) $ as solutions of 
\begin{align}
a_{\rm l}(x^*)(v_{\rm max}-v(x^*))-S(x^*)(1-e^{-a_{\rm l}(x^*)x^*})&=0,\label{eq1-theo1}\\
a_{\rm r}(x^*)v(x^*)-S(x^*) (1-e^{-a_{\rm r}(x^*) (L-x^*)})&=0.\label{eq2-theo1}
\end{align}
The closed-loop  system consisting of  (\ref{nondelay-sys}) with an initial condition $x_0 \in [0,L)$ and the extended  control law   \eqref{left-control}--(\ref{extend-control})  \textcolor{black}{ is GES at  $x=x^*$.}
\label{theoreme-nondelay}
\end{thm}

\textbf{Proof:}  We rewrite the  delay-free plant  (\ref{nondelay-sys})  as 
\begin{align}
\dot e(t)=&   \frac{c(t) \left(U(t)-v(x^*)\right)}{ \left(1-v(x^*)\right)\left(1 +\theta_2(e(t)+x^*)\right)\left(1-U(t)\right)}\nonumber\\  &-\frac{ c(t)\theta_2 e(t)}{1 +\theta_2(e(t)+x^*) },\label{delay-sys-perturb-rewrite}
\end{align}
where, $e(t)= x(t)-x^*$.
The control law (\ref{extend-control}) is written as  
\be
v(x(t),x^*)=  v(e(t)+ x^*,x^*). \label{control-reform}
\ee
 The extended control law   (\ref{extend-control}) is a decreasing function of $x$ and consequently for all $x(t) \in [0,L)$,
\be
\sgn \left (v(e(t)+ x^*,x^*) -v(x^*,x^*)\right)=-\sgn \left(e(t)\right). \label{signe}
\ee
Moreover,  $0\leq v(x,x^*)\leq v_{\rm max}<1$ and $v(x^*,x^*)=v(x^*)$  is  defined as the setpoint open-loop control (\ref{equilibrium}). 
Next, we introduce the following  Lyapunov function
\be
 V
=|e(t)| \label{3}.
\ee
\textcolor{black}{Hence,
\be
\dot V
=\dot e(t) \sgn \left(e(t)\right), \label{4}
\ee
and with the help of  (\ref{signe}), \textcolor{black}{ by choosing $U(t)=v(x(t),x^*)$, (\ref{4}) is written as} 
\be
\dot V=  -\alpha (t)V- \beta(t),\label{dot-V-1}
\ee
where $\alpha (t)>0$ and $\beta(t)\geq 0$  for all $x\in[0,L)$ and  $U\in[0,1)$. The functions $\alpha(t)$ and $\beta(t)$  are   given by  
\textcolor{black}{\begin{align}
\alpha(t) &= \frac{ c(t)\theta_2}{1 +\theta_2(e(t)+ x^*)}\\
\beta(t)&=  \f{1}{\left(1-v(e(t)+ x^*,x^*)\right)}\nonumber\\ & ~~\times \frac{c(t) \left |v(e(t)+ x^*,x^*)-v(x^*)\right|}{\left(1-v(x^*)\right)\left(1 +\theta_2(e(t)+ x^*)\right)},
\end{align}}
where, $ c(t)\geq \theta_1\left(1-\epsilon\right)$, for all $t\geq 0$. Therefore
\be
\dot V
\leq- \f{\theta_1 (1-\epsilon)}{1-\theta_2 L} V.\label{final-lyap-end}
\ee
 From (\ref{final-lyap-end}) the closed-loop system  is exponentially stable at $x^*\in[0,L)$ for all  $x_0 \in[0,L)$. \cqfd }


\section {Predictor Feedback Control for Nonlinear Systems with Time- and State-Dependent Input Delay} \label{sectionpredict} 
\subsection{Predictor feedback design}

We consider the following  nonlinear system with a time- and state-dependent input delay  
\begin{eqnarray} 
\dot x(t) &=& f\left(t, x(t), U(\phi (t)\right) \label {delay-1} \\
\phi (t)  &=& t-D\left(t,x(t)\right), \label {delay-12}
 \end{eqnarray}
where $x \in \mathbb{R}^n$, $U: [t_0-D\left(t_0,x(t_0)\right), \infty) \rightarrow \mathbb{R}$, $t \geq t_0 \geq 0$, $D \in C^1\left(  \mathbb{R}_+\times  \mathbb{R}^n;  \mathbb{R}_+\right)$, and   $ f: R_+ \times   \mathbb{R}^{n} \times   \mathbb{R}\rightarrow \mathbb{R}^n$ is locally Lipschitz with $f\left(t,0,0 \right)=0$ for all $t \geq 0$  and there exists a   class $\mathcal{K}_{\infty}$ function $\hat \alpha $ such that
\begin{eqnarray}\label{f-restrict}
|f\left(t,x,U\right)| \leq \hat \alpha \left( |x|+ |U| \right).
\end{eqnarray}
The predictor feedback control law for system  \eqref {delay-1},  \eqref {delay-12} is
\begin{eqnarray}
U(t)&=&\kappa(\sigma(t), P(t)), \label{control-law} 
\end{eqnarray}
where, for all $t­-D(t,x(t))\leq\theta\leq t$
\begin{align}
P(\theta) &=x(t)+ \int_{t-D(t,x(t))}^{\theta} \frac{ f(\sigma(s),P(s), U(s))}{1-F(\sigma(s),P(s),U(s))} ds, \label{predict-delay-5}\\
\sigma(\theta) &=t+ \int_{ t-D(t,x(t))}^{\theta} \frac{ 1}{1-F(\sigma(s),P(s),U(s))} ds,\label{time-predict}
\end{align}
and
\begin{align}
 F(\sigma(\theta),P(\theta),U(\theta))=&\f{\partial D}{\partial t}\left(\sigma(\theta),P(\theta )\right)\nonumber\\&+\f{\partial D}{\partial x} \left(\sigma(\theta),P(\theta)\right)\nonumber\\&\times f(\sigma(\theta),P(\theta), U(\theta)).\label{F}
 \end{align}

When simulating the predictor feedback controller (\ref {predict-delay-5})--(\ref{control-law}), at each time step the ODE for the system  (\ref{delay-1})--(\ref{delay-12}) must  be solved  (using, for example, a simple Euler scheme) and  the length of the delay  must  be computed (for example as the integer part of $ N(i)=\frac{D( i,x(i))}{\tau}$,  say  $\bar N(i)$, where  $\tau$ is the discretization step). The predictor is then  computed by integrating simultateously the two integral relations (\ref{predict-delay-5}) and (\ref{time-predict})  at each time step,  using a numerical integration scheme. For instance, with the left endpoint rule of integration we get
\begin{eqnarray}\label{num-implement}
P(i)=x(i)+\tau\sum_{k= i-\bar N(i)}^{k=i-1} \frac{f(\sigma(k),P(k),U(k))}{1- F( \sigma(k),P(k),U(k))},\label{num-implement-1}\\
\sigma(i)=i+\tau \sum_{k=i-\bar N(i)}^{k=i-1} \frac{1}{(1- F ( \sigma(k),P(k),U(k))}\label{num-implement-2}.
\end{eqnarray}

The prediction of the state at the time when the current control will have an effect on the state is defined as
\be\label{PP}
P(t)=x\left(t+ D(\sigma(t),P(t))\right),\ee
where the prediction time is defined as
 \be \label{TT} \sigma(t)=t+ D(\sigma(t),P(t)),\ee
which is derived from the inversion of the time variable 
   $t  \to t- D(t,X(t))$ in  $t  \to t+D(\sigma(t),P(t))$, \cite{BEK13,Bekiaris13}.  Differentiating \eqref{PP},  \eqref{TT} and using   \eqref {delay-1} we arrive at
\begin{eqnarray} \label{delay-3}
 \frac{d x(\sigma(t))}{dt}&=& f(\sigma(t), x(\sigma(t)), U(t)) \frac{d\sigma(t)} {dt}, 
 \end{eqnarray}
and
\begin{eqnarray}\label{sigma-deriv}
\dot \sigma(t)&=&\f{1}{1-F(\sigma(t),P(t), U(t))}
 \end{eqnarray}
where $F$ is defined in \eqref{F}.
Finally, the implicit integral relations \eqref{predict-delay-5} and \eqref{time-predict} are derived by integrating  \eqref{delay-3} and \eqref{sigma-deriv} on the delay interval $[\phi(t),\theta]$.

The key point of the predictor feedback  design is the feasibility condition defined as
 \begin{eqnarray} \label{feasibility}
\mathcal{F}_{c}&:&\quad \f{\partial D}{\partial t}\left(\sigma(\theta),P(\theta) \right)\nonumber\\ &&+\f{\partial D}{\partial x} \left(\sigma(\theta),P(\theta)\right) f(\sigma(\theta), P(\theta),U(\theta))< c,
\end{eqnarray}
 for all   $\theta \geq t_0-D(t_0,x(t_0))$ and some  $ c\in (0, 1)$. Condition \eqref{feasibility} guarantees that the feedback control action can reach the plant, namely, the delay rate is bounded by unity, and that the denominator of the predictor    (\ref{predict-delay-5}) and the prediction time (\ref{time-predict}) is positive. We refer the reader to \cite{Bekiaris13}  for details on the predictor feedback control design and analysis for systems with state-dependent input delay.

\subsection{Stability analysis}

\begin{assum}\label {assumption-2}
  There exist   a smooth positive definite function $R$ and  class $\mathcal{K}_{\infty}$ functions $\mu_1$,  $\mu_2$  and   $\mu_3$ such that for the plant $\dot x=f(t,x,w)$, the following hold 
\be
\mu_1(|x|)\leq R(t,x)\leq \mu_2(|x|)\label{FC1}
\ee
\begin{align}
&& \frac{\partial R(t,x)}{\partial t} +\frac{\partial R(t,x)}{\partial x} f(t,x, \omega) \leq R(t, x) +\mu_3 (|\omega|)\label{FC2},  
\end{align}
for all  $(x, \  \omega)^T \in  \mathbb{R}^{n+1}$  and  $ t\geq t_0$.
\end{assum}
Assumption \ref{assumption-2} guarantees that system $\dot{x}=f(t,x,\omega)$ is strongly forward complete with respect to $\omega$.
\begin{assum}\label {assumption-4}
  There exist   a  locally Lipschitz function  $\kappa \in \left([t_0,\infty) \times \mathbb{R}^n ; \mathbb{R} \right)$ and  a function  $\hat \rho \in \mathcal{K}_{\infty}$  such that the plant $\dot x=f\left(t,x(t),\kappa(t,x(t))+\omega(t)\right)$ is input to state stable with respect to $\omega$ and $\kappa$ is uniformly bounded with respect to its first argument,  that is,
\be
|\kappa(t,x)|\leq \hat \rho(|x|) \quad \forall \quad t\geq t_0.
\ee

\end{assum}

\begin{assum}\label {assumption-3}
$D \in C^1(R_+\times R^n ; R_+)$, $\f{\partial D}{\partial t}$ and  $\f{\partial D}{\partial x}$   are locally Lipschitz (to guarantee the uniqueness of solutions), and there exist  class $\mathcal{K}_{\infty}$ functions $\mu_4$, $\mu_5$, $\mu_6$  and  non-negative constants    $c_1$,  $c_2$, $c_3$,  with $c_3<c$, for some $0<c<1$, such that
\begin{eqnarray}
D(t,x(t))&\leq& c_1+\mu_4(|x|),  \label{Feasibility-region-condition1}\\
\left |\f{\partial D}{\partial x}(t,x(t))\right |&\leq 
&c_2+\mu_6(|x|), \label{Feasibility-region-condition3}\\
\left |\f{\partial D}{\partial t}(t,x(t)) \right |&\leq& c_3+\mu_5(|x|). \label{Feasibility-region-condition2}
 \end{eqnarray} 
\end{assum}
The definitions of strong forward completeness and input-to-state stability are those from \cite{Sontag1995}, and \cite{Krstic2010aa}, respectively. 

\begin{thm}\label{theorem-predict}
Consider the closed-loop system  consisting of the plant  (\ref {delay-1}) and  the control law (\ref{predict-delay-5})--(\ref {control-law}). Under Assumptions  \ref {assumption-2}, \ref {assumption-4} and \ref {assumption-3},  there exist a  class $\mathcal{K}$ function $\psi_{\rm RoA} $ and a class $\mathcal{KL}$ function $\beta_s$ such that for all initial conditions for which $U$ is locally Lipschitz on the interval $[t_0 -D(t_0,x(t_0)),t_0)$ and which satisfy
\begin{eqnarray} \label{Feasibility-region-theo}
 \Omega (t_0) <  \psi_{\rm RoA}(c-c_3),
 \end{eqnarray} 
for some $0 < c < 1$, where
\begin{eqnarray} \label{predict-delay-5bis-2}
\Omega(t)=  | x(t)  | +\sup_{ t-D(t,x(t))\leq \theta \leq t }| U(\theta)|,
 \end{eqnarray}
 there exists a unique solution to the closed-loop system with $x$ Lipschitz on $ [t_0,\infty)$, $U$ Lipschitz on $ (t_0,\infty)$, and the following holds
\begin{eqnarray} \label{delay-estimate}
\Omega (t)\leq \beta_s \left(\Omega(t_0),t-t_0\right),
 \end{eqnarray}
for all $ t \geq t_0$. Furthermore, there exists  a positive constant  $\gamma$  such that for all $ t \geq t_0$,
\begin{eqnarray} \label{dot-delay-estimate}
D(t,x(t))&\leq& \gamma\\
\left | \f{d  D(t,x(t))}{dt}\right| &\leq& c\label{dot-delay-estimate-2}
 \end{eqnarray}
\end{thm}

\textbf{Proof of Theorem \ref{theorem-predict}:}
Estimates    (\ref{delay-estimate}), (\ref{dot-delay-estimate}), and (\ref{dot-delay-estimate-2}) follow by directly applying Lemmas 1--8  from  \cite{Bekiaris13} (see the Appendix). Existence and uniqueness of a solution $x$ Lipschitz on $[0,\infty)$ follows from the proof of Theorem 1 in  \cite{Bekiaris13} (page 7). It remains to show that $U$ is Lipschitz on $(t_0,\infty)$.
 Since $U(t)=\kappa\left(\sigma(t), P(t)\right)$  and 
 \begin{align}
\dot P(t) = &\frac{ f\left(\sigma(t),P(t), \kappa(\sigma(t),P(t))\right) }{1- F(\sigma(t),P(t),\kappa\left(\sigma(t), P(t)\right))}, \label{deriv-target-predict}\\
\dot \sigma(t) =& \frac{ 1}{1-F(\sigma(t),P(t),\kappa\left(\sigma(t), P(t)\right)) },\label{ deriv-sigma-target-predict}\\
F(\sigma(t),P(t))=&\f{\partial D}{\partial x} D\left(\sigma(t),P(t)\right)\nonumber\\ &\times f\left(\sigma(t),P(t), \kappa\left(\sigma(t), P(t)\right))\right) \nonumber\\ &+ \f{\partial D}{\partial t}\left(\sigma(t),P(t) \right),
 \end{align}
for $t\geq t_0$, the Lipschitzness of  $\f{\partial D}{\partial t}$, $\f{\partial D}{\partial x}$, $ \kappa$ and $f$,  and  (\ref{feasibility}) ensure that the right hand-side of \eqref{deriv-target-predict} and \eqref{ deriv-sigma-target-predict} are Lipschitz and consequently $(P ,\sigma)\in( C^1(t_0, \infty)\times C^1(t_0, \infty))$. From the Lipschitzness of  $\kappa $, it follows that $U$ is  Lipschitz.\cqfd


\section{Application to the Extrusion Process Model} \label{application} From now, we recall the predictor feedback   (\ref{predict-delay-5})--(\ref {control-law}) for the compensation of the time- and   state-dependent input delay  in system (\ref{delay-sys-perturb}) that we rewrite  formally as 
\begin{eqnarray}\label{predicteur}
\dot x(t) &=& f \left(t, x(t), U\left(t- D(t,x(t))\right) \right)\label{delay},\\
f \left(t, x(t), U(t) \right) &=& - c(t) \Gamma(x(t),U(t)),\label{delay-sys-perturb-bis}
\end{eqnarray}
\textcolor{black}{where  $D(t,x(t))$ and  $c(t)$ are defined in (\ref{delay-time}) and \eqref{delay-time-fluc}, respectively, and}
\begin{align}
&\Gamma(x(t),U(t))=\left[\frac{ \theta_2  x(t)}{\left(1 +\theta_2x(t)\right) \left(1-U\left(t\right)\right)}-\frac{ U \left(t)\right)}{ \left(1-U\left(t\right)\right)}\right],\nonumber\\
& U(t) \in [0,1). \label{nondelay-sys-grad}
\end{align}
The predictive feedback  controller based on the piecewise exponential feedback law (\ref{extend-control}) is given by  
\begin{align}
& U(t)=
v\left( P(t)\right),\label{predictive-controller}\\
&P(\theta) =x(t)+ \int_{t-D(t,x(t))}^{\theta} \frac{ f(\sigma(s),P(s), U(s))}{1-F(\sigma(s),P(s),U(s))} ds \label{predict-delay-extru}\\
&\sigma(\theta) =t+ \int_{ t-D(t,x(t))}^{\theta} \frac{ 1}{1-F(\sigma(s),P(s),U(s))} ds, \label{control-law-extru}\end{align}
for all $t-D(t,x(t))\leq \theta \leq t$. \textcolor{black}{The function $F$ defined in  (\ref{F}) for the system  (\ref{delay})--(\ref{delay-sys-perturb-bis})   is computed with the help of  (\ref{delay-time})--(\ref{delay-time-fluc}) as}
\begin{align}
 F(\sigma(t),P(t),U(t))=& \f{\theta_1 \epsilon \omega \sin(\omega \sigma(t))(L-P(t))}{  c^2(t)} \nonumber \\ &+
  \Gamma(P(t),U(t)), \label{extru-Feas}
\end{align}
where $\Gamma(P(t),U(t))$ is defined in \eqref{nondelay-sys-grad}.

The parameters $a_{\rm l}(x^*)$ and $a_{\rm r}(x^*)$  of the feedback control law \eqref{predictive-controller} are the  solutions of  (\ref{left-trans-eq}) and  (\ref{right-trans-eq}) for an assigned slope function value at the set point that satisfies (\ref{condition-domain}).
$P(t)=x\left(t+ D(\sigma(t), P(t))\right)$ is the prediction of the state at the time when the current control will have an effect on the state.  Recall that  the implicit  integral relation  (\ref{predict-delay-extru}) is derived from the inversion of the time  variable $t  \to t- D(t,x(t))$ in  $t  \to t+D(\sigma(t),P(t))$  with the prediction time defined as $\sigma(t)=t+ D(\sigma(t), P(t))$. The key point of the design is the feasibility condition $\mathcal{F}_c$ defined in  (\ref{feasibility}),  which ensures that the control action can reach the plant, namely,  the delay rate is bounded by unity. The a priori satisfaction of (\ref{feasibility}) depends on the magnitude $\epsilon$ and the angular frequency $\omega$  of the periodic instability,  and on  the design parameters of the extruder.

\begin{thm}\label{theorem-predict-gas}
 For any setpoint $x^* \in (0,L)$ and for any  choosen  setpoint slope  $S(x^*)\in \mathbb{R}$    satisfying $S(x^*)\geq S_{\rm min}(x^*)$, where $S_{\rm min}(x^*)$  is given by (\ref{condition-domain}) and any initial condition $x_0 \in [0,L)$ and,
\begin{eqnarray}
\left \{U_0(\theta)\, | \,U_0(\theta)\in [0,1), \,\,\,\textrm{for all}\,\, \theta \in [-D(t_0,x_0),0) \right\},
\end{eqnarray}
taking the control gains $a_ {\rm l}(x^*)$ and $a_ {\rm r}(x^*) $ as solutions of (\ref{eq1-theo1})  and  (\ref{eq2-theo1}), respectively,  the closed-loop  system consisting of   the plant (\ref{delay})--(\ref{nondelay-sys-grad})  with state $x(t)$, together with the control law (\ref{predictive-controller})--(\ref{extru-Feas}), \eqref{extend-control} with actuator state $U(t+\theta)$, $\theta \in [-D(t,x(t)),0)$,  is \textcolor{black}{GAS at  $x=x^*$,  $U=v(x^*)$}
 if the parameters of the extruder model and the perturbation satisfy, 
\textcolor{black}{\begin{eqnarray}
0\leq \f{\epsilon  \omega}{(1-\epsilon)^2}&<& \frac{\theta_1\theta_2  }{  (1 +\theta_2 L)^2}, \label{eq-feas-4}
\end{eqnarray}}
\textcolor{black}{or, 
\begin{eqnarray}
\frac{\theta_1 \theta_2  }{ (1 +\theta_2 L)^2} <\f{\epsilon \omega }{(1-\epsilon)^2}<\f{\theta_1}{L},  \quad \textrm{and} \quad \theta_2<\f{1}{L}, \label{feas--decreas-theo-bis}
\end{eqnarray}
 or, 
\be
 \frac{\theta_1 \theta_2  }{ (1 +\theta_2 L)^2}<\f{\epsilon \omega }{ (1-\epsilon)^2}< \frac{4\theta_1 \theta_2  }{ (1 +\theta_2 L)^2}, \quad \textrm{and} \quad  \theta_2>\f{1}{L},\label{condi-feas-last-theo}
\ee}
where, $\theta_1$ and $\theta_2$ are defined in     \eqref{theta-1-theo} and \eqref{theta-2-theo}, respectively.
\end{thm}

\textbf{Proof:} The proof of  Theorem \ref{theorem-predict} is based on  the Lyapunov-like condition   (\ref{feasibility}) that must be satisfied   $a$  $priori$ to guarantee the GAS property for any  given $x^*\in(0,L)$. In the following, we compute the function (\ref{F}) for   the time- and  state-dependent input  delay model   of the extruder (\ref{delay})--(\ref{nondelay-sys-grad}), (\ref{delay-time}), \eqref{delay-time-fluc}  in order to establish that the feasibility region as it is defined by (\ref{feasibility}) is the entire physical domain, namely, $x\in(0,L)$ and $U\in[0,1)$. It holds that
\begin{align} \label{Time-deriv}
\f{\partial D}{\partial t}(t,x)&= \f{\theta_1 \epsilon \omega \sin(\omega t)}{ c^2(t)}(L-x),\\
\f{\partial D}{\partial x}(t,x) f(t,x,U)&=\Gamma(x,U),\label{grad}
\end{align}
where $\Gamma(x,U)$ is given by \eqref{nondelay-sys-grad}. Since, $0<\epsilon<1$,  for all  $x\in(0,L)$ we get that
$
\min c(t)=\theta_1 (1-\epsilon),
$
and hence,
\begin{eqnarray} \label{Time-deriv-bound}
\f{\partial D}{\partial t}(t,x) \leq \f{\epsilon \omega(L-x)}{ \theta_1 (1-\epsilon)^2}.
\end{eqnarray}
The gradient of \eqref{nondelay-sys-grad}  with respect to the  input $U$ satisfies
\be
\nabla_U \Gamma (x,U)
= - \frac{\theta_2 }{\left(1+\theta_2x\right)\left(1-U\right)^2}.\label{phi-grad}
\ee
It follows that  (\ref{nondelay-sys-grad}) is a  strictly decreasing function of $U$ which belongs to  $[0,v_{max}]$,  for all $x(t) \in [0, \infty)$ and 
\begin{eqnarray}
\sup_{U \in [0,v_{max}]}{\Gamma(x,U)} &=& \frac{\theta_2  x}{(1 +\theta_2 x)}. \label{nondelay-sys-grad-1}
\end{eqnarray}
 The delay rate $\f{d D}{d t}$ is uniformely bounded by unity,  namely,  the feasibility condition (\ref{feasibility})   is satisfied if and only if
\begin{eqnarray} \label{Total-deriv}
\f{\partial D}{\partial t}(t,x)+
\f{\partial D}{\partial x}(t,x) f(t,x,U)<1.
\end{eqnarray}
for all  $t\geq 0$,  $x\in(0,L)$ and $U\in[0,1)$.  By  (\ref{Time-deriv-bound}) and  (\ref{nondelay-sys-grad-1}), it follows that (\ref{Total-deriv}) is satisfied   for all $x \in (0,L)$ and  $U \in [0,v_{max}]$    
\begin{enumerate}
\item if  \begin{eqnarray}
\Lambda(x)= \f{\epsilon \omega (L-x)}{\theta_1 (1-\epsilon)^2}+\frac{\theta_2  x}{(1 +\theta_2 x)}, \label{eq-feas-1}
\end{eqnarray}
is a strictly increasing function since then  its  maximum  over the domain $(0,L)$ satisfies   
\textcolor{black}{ \begin{eqnarray}
\Lambda(L) =\frac{\theta_2  L}{(1 +\theta_2 L)}<1.  \label{nondelay-sys-grad-L}
\end{eqnarray}}
The derivative of the function (\ref{eq-feas-1}) is written as
\begin{eqnarray}
\Lambda^{'} (x)=\frac{\theta_2  }{(1 +\theta_2 x)^2}-\f{\epsilon \omega }{ \theta_1(1-\epsilon)^2}, \label{eq-feas-2}
\end{eqnarray}
and hence, since $x\in(0,L)$,  (\ref{Total-deriv})  is guaranted if \eqref{eq-feas-4} holds.

\item if   \eqref{eq-feas-1} is   a decreasing function we deduce from \eqref{eq-feas-2} that
\begin{eqnarray}
\f{\epsilon \omega }{ (1-\epsilon)^2}> \theta_1 \theta_2,  \label{eq-feas-3-decrease}
\end{eqnarray}
and
\begin{eqnarray}
\sup_{x \in [0, L]}{\Lambda(x)} &=& \f{\epsilon \omega L}{\theta_1 (1-\epsilon)^2}. \label{nondelay-sys-grad--decreas}
\end{eqnarray}
Finally, the feasibility condition is satisfied if 
\begin{eqnarray}
\theta_1 \theta_2 <\f{\epsilon \omega }{(1-\epsilon)^2}<\f{\theta_1}{L}. \label{feas--decreas-}
\end{eqnarray}
One should notice that \eqref{feas--decreas-} necessarily restricts $\theta_2$ to satisfy \be\theta_2<\f{1}{L}\label{restrict-bis}.\ee 

\item if   \eqref{eq-feas-1} is an increasing function of $x$ on  the interval $[0,x_1]$  and a decreasing function on $[x_1, L]$ such that $\Lambda(x_1)<1$. The   maximum value of $\Lambda$ is attained at $x_1$ satisfying
\be
x_1= (1-\epsilon) \sqrt{\f{\theta_1}{\theta_2 \epsilon \omega}}-\f{1}{\theta_2}.
\ee
\textcolor{black}{It becomes clear that \eqref{eq-feas-1}  admits a unique maximum  satisfying the feasibility condition if 
\be\label{condi-1}
0< (1-\epsilon) \sqrt{\f{\theta_1}{\theta_2 \epsilon \omega}}-\f{1}{\theta_2} < L,
\ee
and                                                                
\be\label{condi-2} \Lambda(x_1)<1.\ee
The relation  \eqref{condi-1} is equivalent to 
\be\label{condi-1-bis}
 \frac{\theta_1 \theta_2  }{ (1 +\theta_2 L)^2}<\f{\epsilon \omega }{ (1-\epsilon)^2}<\theta_1 \theta_2.  
\ee
Using  \eqref{eq-feas-1}, the inequality \eqref{condi-2} leads to the following relation
\be
\f{\epsilon \omega }{ (1-\epsilon)^2}< \frac{4\theta_1 \theta_2  }{ (1 +\theta_2 L)^2}\label{condi-2-bis}.\ee
For satisfying \eqref{condi-1-bis} and \eqref{condi-2-bis} we need to either  impose \eqref{condi-feas-last-theo} or the condition \eqref{condi-1-bis} with $\theta_2<\f{1}{L}$ which can be combined with  \eqref{feas--decreas-} and \eqref{restrict-bis} in order   to derive \eqref{feas--decreas-theo-bis}. \cqfd}
\end{enumerate}

\begin{figure}[!h]
\begin{center}
\includegraphics[height=140\unitlength]{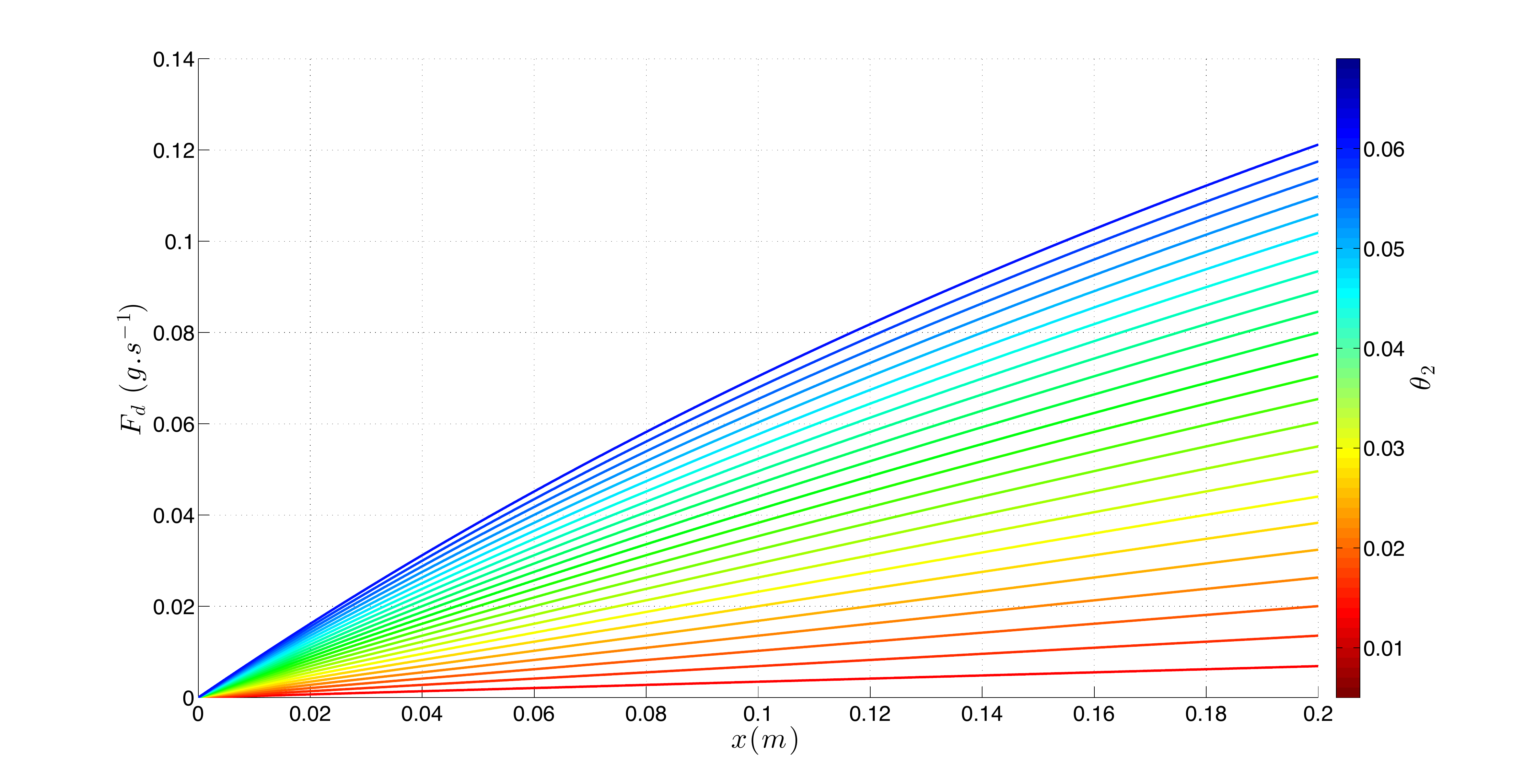}
\caption{\label{Fd} The nozzle flow rate  is an increasing  function of $\theta_2$.}
\end{center}
\end{figure}

 \begin{remark}
\textup{  For given values for the parameters of the extruder, namely, $\theta_1$, $\theta_2$, and
 $L$, the condition (\ref{eq-feas-4})  is always satisfied  if   $ \epsilon$ or $\omega$  are sufficiently small. An increase of the magnitude of   $\epsilon$  causes a decrease of the allowed  $\omega$   and vice versa, as it is evident from (\ref{eq-feas-4})--\eqref{condi-feas-last-theo}. For 
 given  $\theta_1$, and $L$, the maximum bound of the perturbation parameters, namely, $\epsilon$ or $\omega$  is   expressed in   (\ref{eq-feas-4}) as
\begin{eqnarray}\label{max}
\sup_{\theta_2 \in \mathbb{R}} \left \{\frac{\theta_1 \theta_2  }{  (1 +\theta_2 L)^2}\right \}=\frac{\theta_1  }{4 L }, \quad \theta_2=\f{1}{L}
		\end {eqnarray}
Larger variations of $\epsilon$ and $\omega$ are possible, especially in the case in which $\theta_2$ is small, as it is evident from  \eqref{feas--decreas-theo-bis} and  \eqref{condi-feas-last-theo}. However, for very large $\theta_2$, one can conclude from  (\ref{eq-feas-4}) and \eqref{condi-feas-last-theo}  that the allowable size of $\epsilon$ and $\omega$ is restricted.
Moreover, from  (\ref{eq-feas-4})--\eqref{condi-feas-last-theo} one can conclude that the size of the allowable fluctuations of the transport speed in $\epsilon$ and $\omega$ is proportional to $\theta_1$ and inversely proportional to the extruder length $L$. \\\\
In physical terms, conditions (\ref{eq-feas-4})--\eqref{condi-feas-last-theo}  are mainly a correlation between  the pressure   and    the ``rotation" flow, namely,   $\theta_2$ defined in (\ref{theta-2-theo}) and  $ \theta_1$ defined in (\ref{theta-1-theo}), respectively. We recall the expression of the net flow rate defined in \eqref{pressure-conti} which is    an increasing function of  $\theta_2$  as it is shown in  Fig. \ref{Fd}. Therefore, changes in $\theta_2$, by manipulating $K_d$, $B$, or $\rho_0$, the nozzle conductance, the screw resistance, and the melt density, respectively,  affect the output flow rate $F_d$. For example, an increase in $\theta_2$  by increasing the nozzle conductance $K_d$, leads to an increase in the outflow rate. Note that $K_d$,  which defines the nozzle opening, is directly related to the printing resolution, namely, the accuracy of the printing process. A large nozzle opening leads to an extrusion of a filament with a large diameter and consequently deteriorates the printer precision. Moreover, from   \eqref{feas--decreas-theo-bis} and  (\ref{theta-2-theo}), it can be also seen that the  ``robustness" of the controller depends on  the material   thickness, namely, the mass density $\rho_0$: a thicker material is  less sensitive to large fluctuations of the \textcolor{black}{transport speed}  under the  predictor feedback control law.   The parameter  $B$  in the expression of $\theta_2$  in \eqref{theta-2-theo} is given by
\be\label{B}
B=\f{WH^3}{12},
\ee
where $H$ is the approximate depth of screw channel  from  the screw thread root to the  barrel internal surface and $W$ is the width of screw channel. Consequently, changes in  $\theta_2$ due to the changes in $B$ affect also   the parameter $\theta_1$,  since  the screw pitch value $\xi$ also depends directly on $W$.\\\\
Relations (\ref{eq-feas-4})--\eqref{condi-feas-last-theo}  show that  an increase in $\theta_1$, namely,  an increase of the material convection speed,  by enabling a large screw pitch $\xi$ or a high screw speed $N_0$,  improves the ``robustness" of the controller in some way and allows for a system that supports broader changes of the convection velocity in frequency and amplitude. Note that  a sharp increase in the rotational  screw speed $N_0$   results  in material overload and clogging  problems and has a major effect on the residence time  that  is the critical  time during which the material should be heated to have good properties before being evicted through the nozzle. Particulary,  the  extruded filament homogeineity is directly related to the residence time and to the  process of solidification after layers deposition  in 3D printers.  In addition, an increase in $\theta_1$ in the screw speed $N_0$,  increases   the thermal  energy in the extruder chamber  due to the material shearing and decreases the viscosity of the melt.  In that case, a rapid feeding of the extruder with  granular material by applying  a more agressive ``Bang-Bang" control action   absorbs  the excess heat in the system.  Maintaining  a reasonnable temperature inside the barrel  is essential because an excessive overheating of the system burns the polymer or produces poor extrusion. Generally, the conventional extrusion processes are equipped with a cooling system to compensate for the heat generated by the mechanical shearing effect that is proportional to the screw speed.\\\\
 In general, the nozzle and the screw designs are directly related to the predictor feedback control design and for achieving high performances for the closed-loop system the scale of the extruder should be neatly chosen.  For instance, the  agressiveness  of the controller is  influenced by the choice of $\theta_2$ since the minimum value of the slope at the set point $S_{\rm min}(x^*)$ defined in \eqref{condition-domain}  depends on   this parameter.
Moreover, the  entire process operates  with an   extruder head  that moves very fast to print filament lines  layer upon layer on a moving platform.  A  sufficiently light extruder head with small nozzle opening $K_d$ and a small length $L$   that operates    at a sufficiently high screw speed $N_0$  is needed to ensure a  high rate of extrusion with  a high precision.}
\end{remark}
   
\section{Simulations}\label{simulation}
\subsection{Time- and state-dependent input delay compensation}  The setpoint is chosen as   $x^*=0.16 \ \rm{m}$   which corresponds to a desired nozzle output flow rate as indicated by the equation \eqref{pressure-conti}.  The initial   position  of the moving interface is set to $x_0=0.1\ \rm{m}$, the  total length of the extruder is  $ L=0.2 \ \rm{m}$ and the system is supposed to settle  at     $x_0=0.1\rm{m}$  at the initial time. The value of the slope function is set to   $S(x^*)= S_{\rm min}(x^*)+30$. The simulations show the dynamics of the input  filling ratio $U(t)$, the interface position $x(t)$, the predictor state $P(t)$ and the delay function $D(t,x(t))$. Different cases including   the open-loop dynamics, both uncompensated and compensated delay control laws are simulated  for $\{\epsilon=0.1, \, \omega=3.5 \   \rm{rad/s}\}$ and $\{\epsilon=0.4, \, \omega=0.4\  \rm{rad/s}\}$.  It is clear that, the uncompensated input leads to a limit cycle and the compensated closed-loop control  allows faster convergence than  the open loop control. Also, as it is shown in Figure \ref{Feasfig}, the feasibility condition is satisfied in both  presented simulation results. 
\begin{figure}[!h]
\setlength{\unitlength}{1.1pt}
\begin{minipage}[t]{1\textwidth}
\includegraphics[height=125\unitlength]{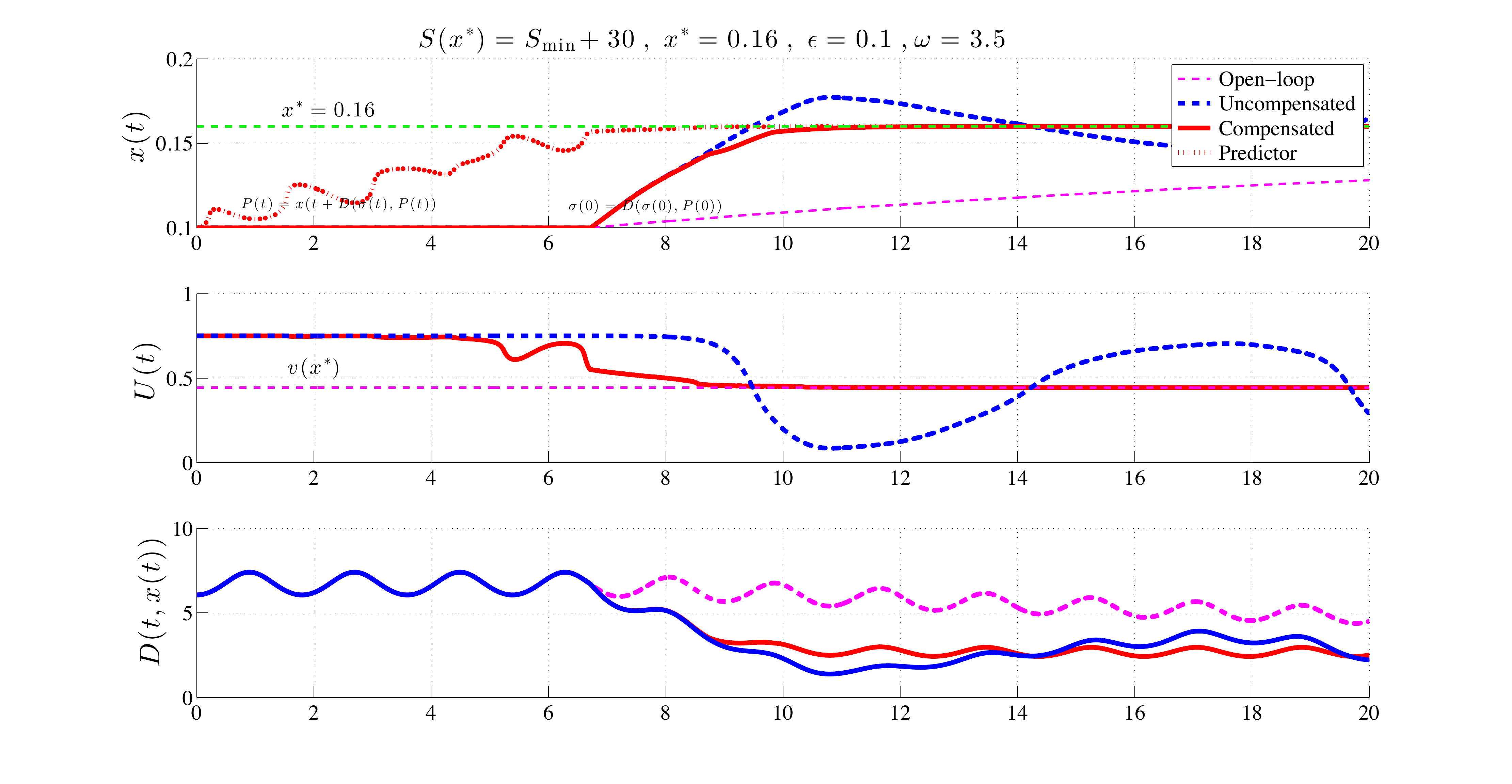}
\caption{\label{Time1} Compensation of the time- and state-dependent input delay-(a).}
\end{minipage}
\end{figure}

\begin{figure}[!h]
\setlength{\unitlength}{1.1pt}
\begin{minipage}[t]{1\textwidth}
\includegraphics[height=125\unitlength]{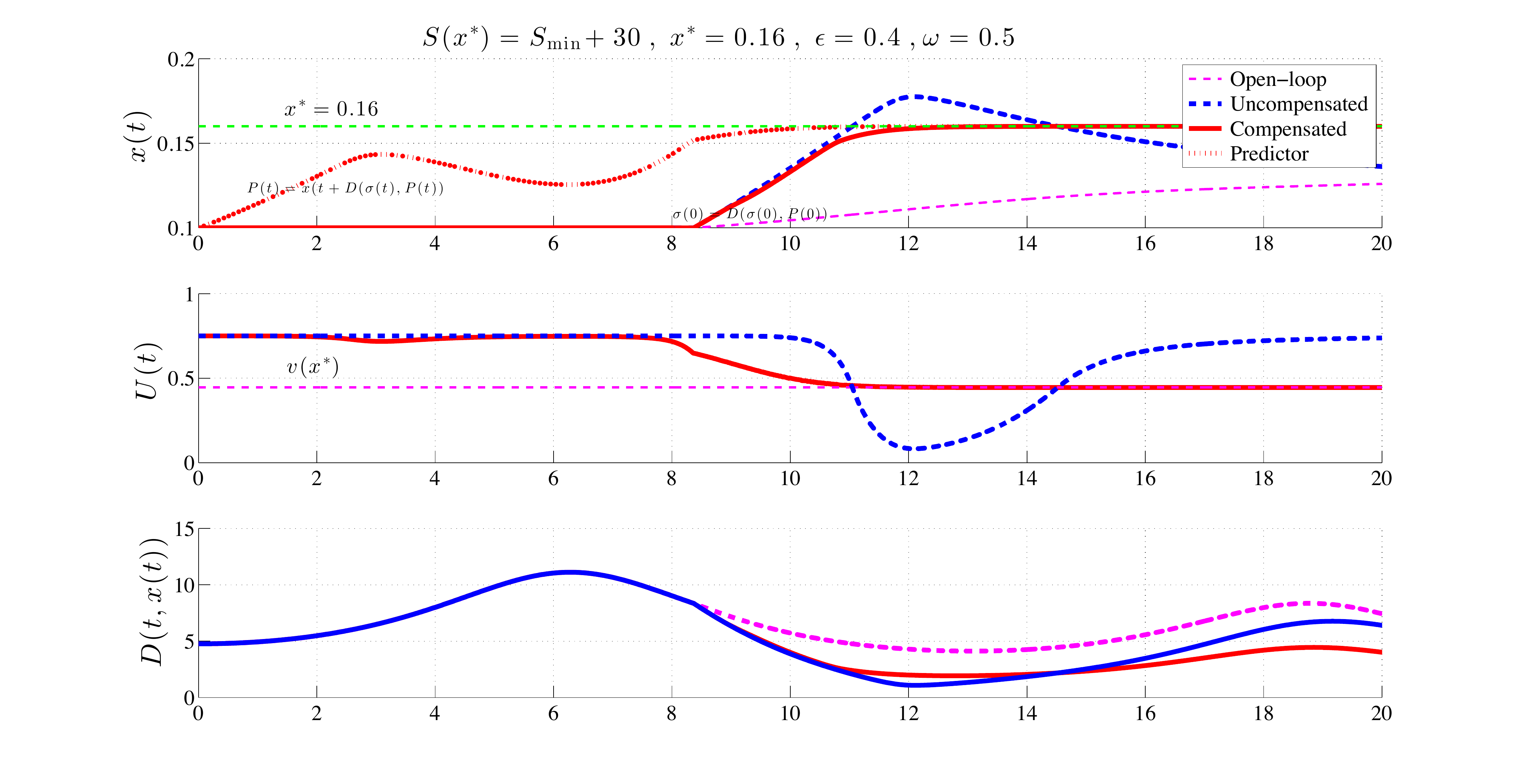}
\caption{\label{Time2}  Compensation of the time- and state-dependent  input delay-(b). }
\end{minipage}
\end{figure}

\subsection{The state-dependent input delay compensator for the model with constant viscosity} The following  simulation results show the stabilization  of the model described by   (\ref{delay-sys}) with  the state-dependent input delay predictor feedback  law  \cite{Bekiaris13}. Defining (\ref{delay-sys})  as 
\begin{eqnarray}\label{sys-normal}
\dot x(t)&=& f\left(x,U\left(t-D_{\rm s}(x(t) \right) \right),\\
D_{\rm s}(x(t))&=& \f{L-x(t)}{\theta_1},
  \end{eqnarray}
where
\begin{eqnarray}\label{sys-normal-1}
f(x(t),U(t))&=&- \theta_1 \Gamma(x(t),U(t)),
  \end{eqnarray}
$\theta_1$ is  the nominal  \textcolor{black}{ transport velocity of the material} defined in  (\ref{theta-1-theo}), and  the function $\Gamma(x(t),U(t))$ is given by  \eqref{nondelay-sys-grad}.  The predictor feedback controller is written as
\begin{align}
&U(t)=
v\left( P_{\rm s}(t)\right),\label{predictive-controller-x}\\
&P_{\rm s}(t) =x(t)+ \int_{t-D_{\rm s}(x(t))}^{t} \frac{ f (P_{\rm s}(\mu), U(\mu))}{1-F_s(P_{\rm s}(\mu), U(\mu))} d\mu \label{state-predict-delay-extru-x}.  \end{align}
where for all $t-D_{\rm s}(x(t))\leq\mu\leq t$
\begin{align}\label{predictive-controller-x-bis}
F_s(P_{\rm s}(\mu), U(\mu))= \f{\partial D_{\rm s}}{\partial x}(P_{\rm s}(\mu))f(P_{\rm s}(\mu), U(\mu)). \end{align}
By specializing Theorem \ref{theorem-predict-gas} to the case $\epsilon=0$ it can be shown that the predictor feedback law   (\ref{predictive-controller-x})-- (\ref{predictive-controller-x-bis}) renders system \eqref{sys-normal}--\eqref{sys-normal-1} GAS (in the physical domain) at any given equilibrium $x^*$.

 \begin{figure}[!ht]
\setlength{\unitlength}{1.1pt}
\begin{minipage}[t]{1\textwidth}
\includegraphics[height=125\unitlength]{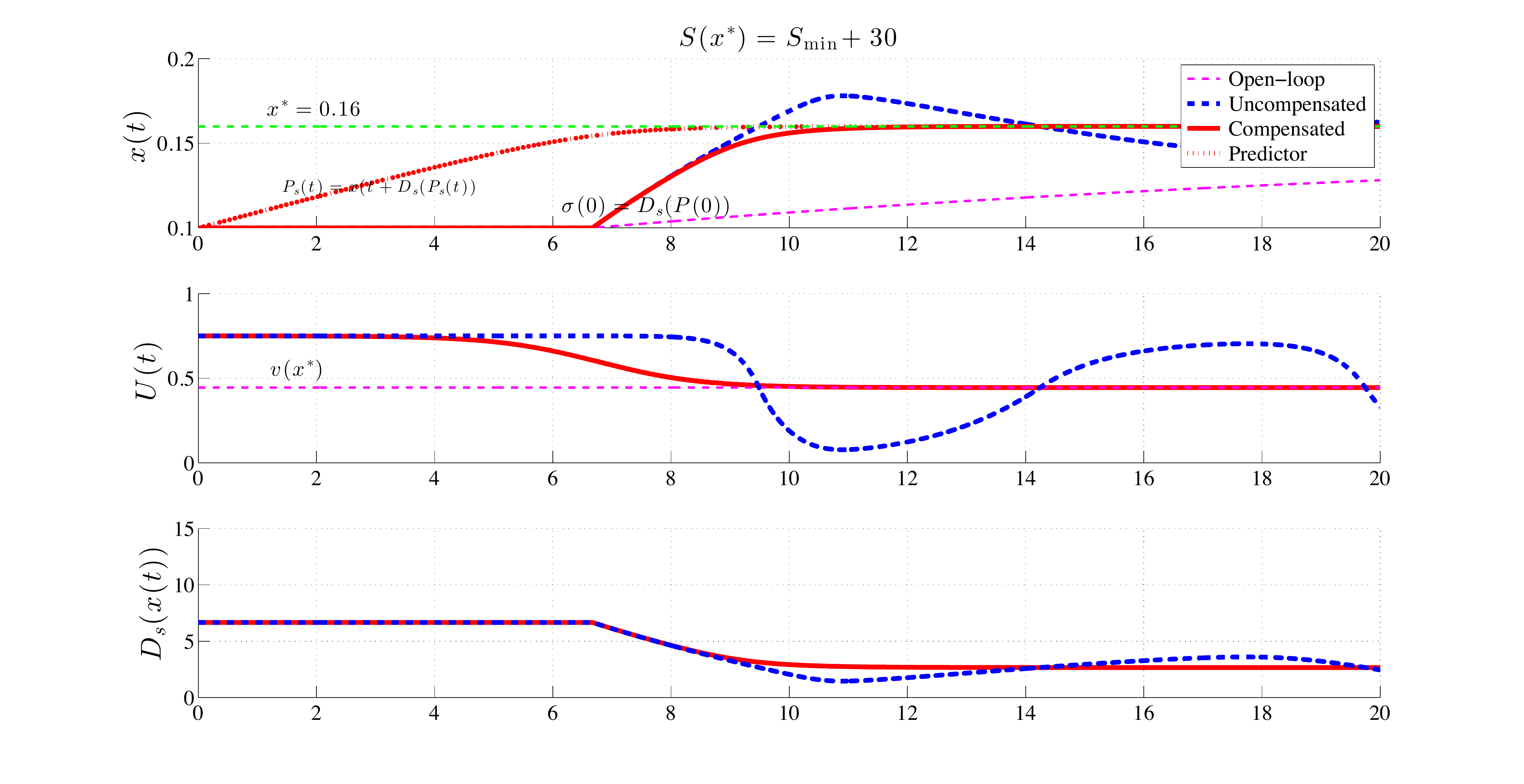}
\caption{\label{State-D(x)} Compensation of the state-dependent input delay.  }
\end{minipage}
\end{figure}

\subsection{Control with a state-dependent input delay compensator} We deal with the case in which the time variations of the transport speed are unknown and  consider  the closed-loop system consisting of   \textcolor{black}{ the plant \eqref{delay-sys-perturb}  with an actual delay $D(t,x(t))$,  given in  (\ref{delay-time}), \eqref{delay-time-fluc}}, together  with a  state- dependent input delay predictor feedback    defined by   
\textcolor{black}{\begin{align}
&U(t)=
 v\left(\hat {P}(t)\right),\label{predictive-controller-1}\\
& \hat{P}(t) =x(t)+ \int_{t-\hat{D}(x(t))}^{t} \frac{ \hat f(\hat{P}(s), U(s))}{1-\hat F(\hat{P}(s),U(s))} ds \label{state-predict-delay-extru-2}, \end{align}
where,
\begin{align}
\hat  F(\hat{P}(s),U(s))= \f{\partial \hat{D}}{\partial x}(\hat{P}(s)) \hat f(\hat{P}(s), U(s)). \end{align}
With an estimated delay function defined as
\begin{eqnarray} \label{hat-delay}
\hat D(x(t))&=& \f{L-x(t)}{\xi N_0},
  \end{eqnarray}
and the nominal vector field,
\begin{align}
\hat f\left(x(t),U(t)\right)=&  \theta_1 \left[-\frac{ \theta_2  x(t)}{\left(1 +\theta_2x(t)\right) \left(1-U\left(t\right)\right)} \right. \nonumber \\&+ \left. \frac{ U \left(t\right)}{ 1-U\left(t\right)}\right], \quad U(t) \in [0,1).\label{delay-sys-estimate}
\end{align}
More precisely, both the predictor  state, $\hat{P}(t)$  and the delay function,    $\hat D(x(t))$,  are estimates of the actual   prediction state, namely, $P(t)$,  and  delay function, namely, $D(t,x(t))$,  that are described by   (\ref{predict-delay-extru}) and  (\ref{delay-time}), \eqref{delay-time-fluc}, respectively.  
For implementing the controller (\ref{predictive-controller-1}), (\ref{state-predict-delay-extru-2}),  the  ``actual" feasibility condition, defined in  (\ref{eq-feas-4}),  \eqref{feas--decreas-theo-bis}, and \eqref{condi-feas-last-theo}  have to  hold, in order to guarantee that the controller actually ``kicks in". In addition, we assume that the following condition, which guarantees that the denominator in (\ref{state-predict-delay-extru-2})  remains always positive (and hence, the controller remains bounded) is satisfied}
\begin{eqnarray}
&& \frac{\theta_2 \hat{P}(\theta)}{\left(1+\theta_2\hat P(\theta)\right) \left(1-U(\theta)\right)}-\frac{ U(\theta)}{ \left(1-U(\theta)\right)}<1,\nonumber\\
&&  \theta\in[t-\hat{D}(x(t)),t].\label{feas-x-estim}
\end{eqnarray}
Note that with strictly positive physical parameters $B$, $\rho_0$, and $K_d$, and for $U\in[0,v_{max}]$, relation (\ref{feas-x-estim}) is satisfied whenever $\hat{P}(\theta)\in [0,\infty)$, for all $\theta\in[t-\hat{D}(x(t)),t]$.

\begin{figure}[!ht]
\setlength{\unitlength}{1.1pt}
\begin{minipage}[t]{1\textwidth}
\includegraphics[height=125\unitlength]{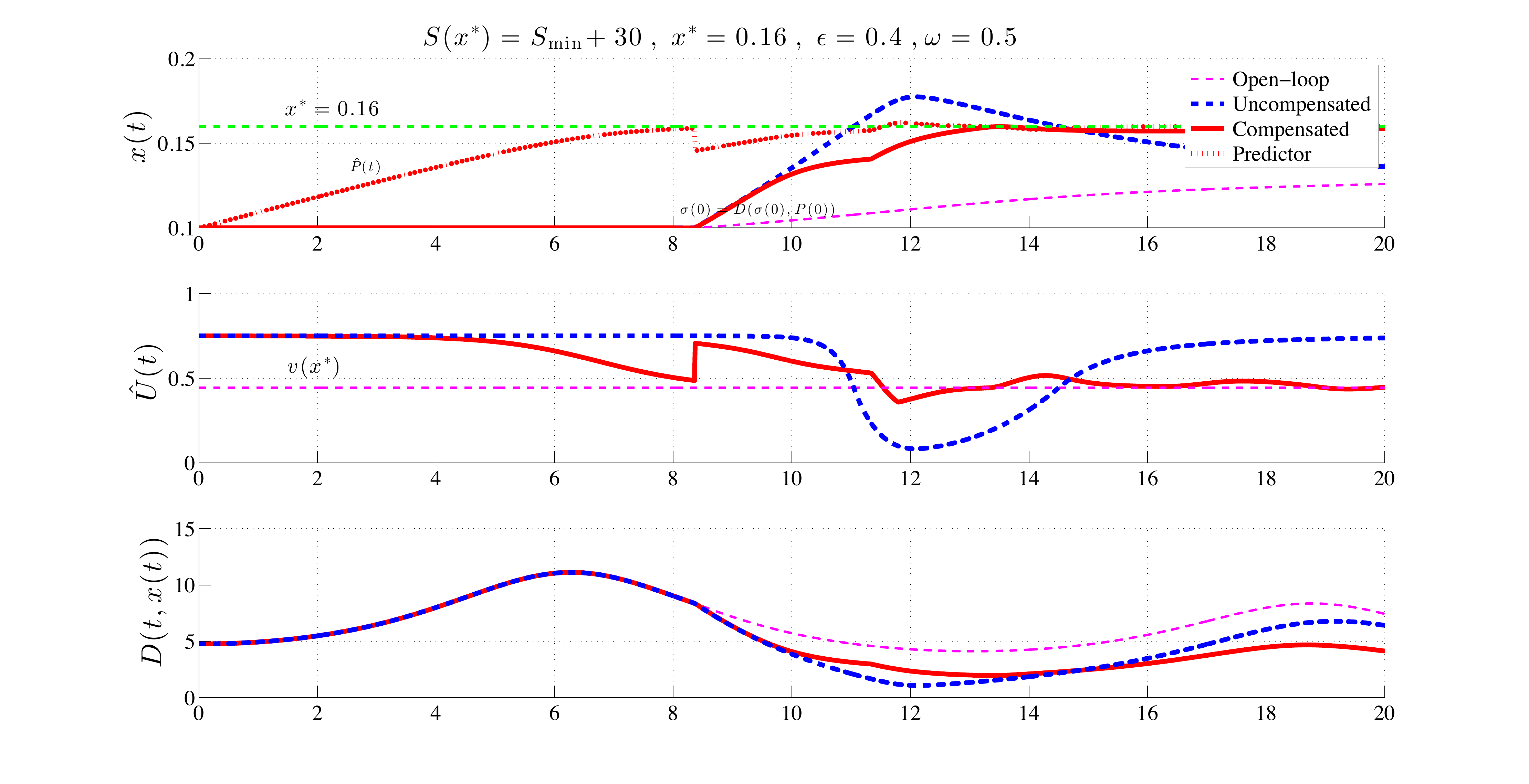}
\caption{\label{State} Robustness of the state-dependent input delay compensator.}
\end{minipage}
\end{figure}

\begin{figure}[!ht]
\setlength{\unitlength}{1.1pt}
\begin{minipage}[t]{1\textwidth}
\includegraphics[height=125\unitlength]{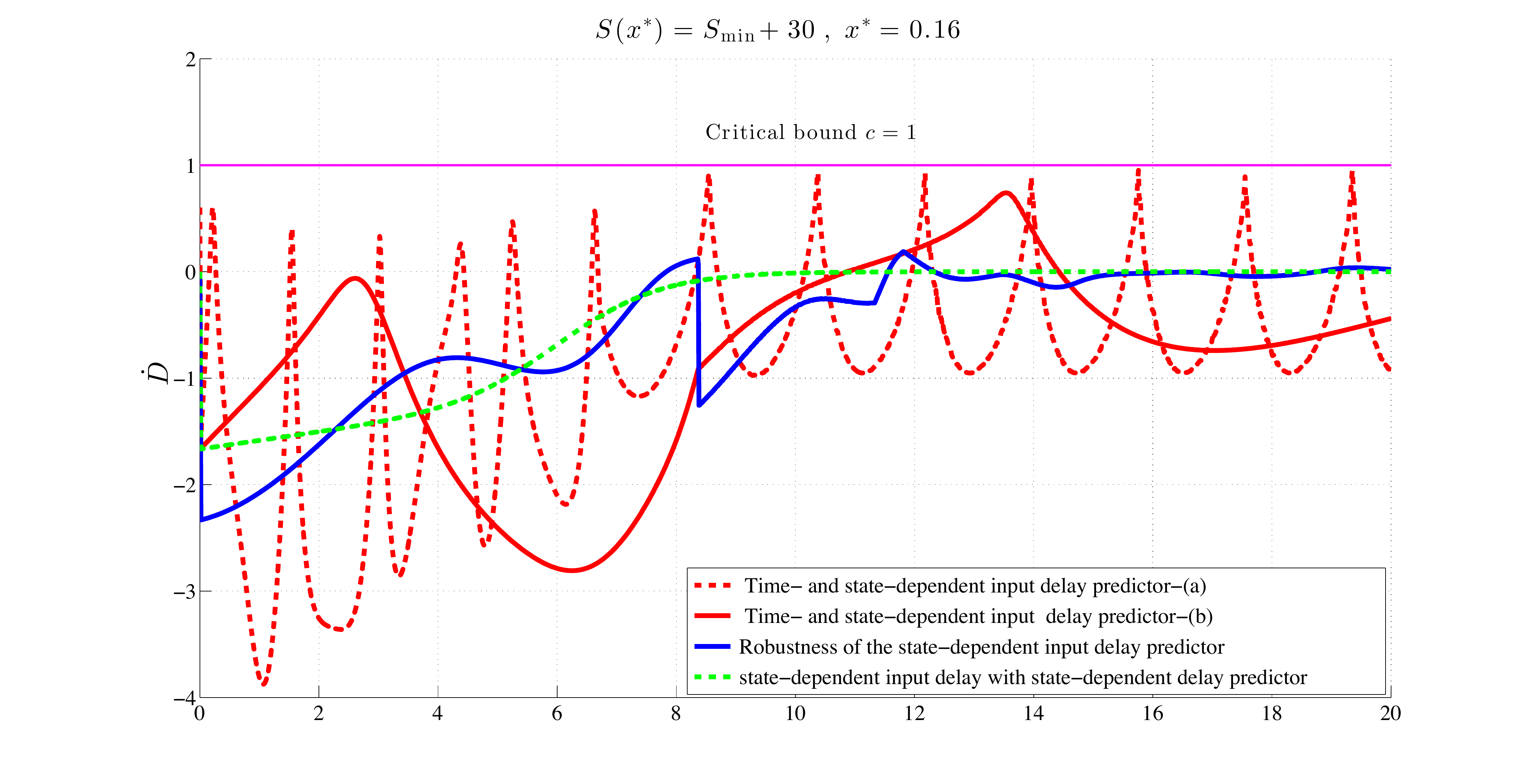}
\caption{\label{Feasfig} Feasibility condition: the delay rate is bounded by unity.  }
\end{minipage}
\end{figure}

The simulation results in Fig. \ref{State} illustrate that the state-dependent input delay compensator can handle small time-varying uncertainties on   \textcolor{black}{the vector field \eqref{delay-sys-estimate}  and the delay function \eqref{hat-delay}, as described by \eqref{delay-sys-perturb-bis}--\eqref{nondelay-sys-grad} and  (\ref{delay-time})--\eqref{delay-time-fluc}, respectively.} An increase in the necessary control effort    to drive the  system to the setpoint is also denoted  and the   rate of convergence decreases  compared to  the time- and state-dependent predictor feedback (\ref{predictive-controller})--(\ref{control-law-extru}) shown in   Fig. \ref{Time2}.

\section{Conclusions}

This paper is devoted to the stabilization of a screw-extrusion process. A coupled PDE-ODE model is used to derive a state-dependent input delay system describing the melt convection in the extruder chamber for an isothermal case. The extension of the aforementioned model to a non-isothermal case is proposed introducing a periodic time-dependent function in the state-dependent input delay function. Next, we design a predictor feedback controller to compensate the state- and time-dependent input delay and establish the GAS of any setpoint with respect to the physical domain under physical and design restrictions. The delay compensator is constructed with a nominal Bang-Bang-like controller that ensures the GES of the delay-free plant.

It is clear that the delay function model for the non-isothermal extrusion process should depend on the rheological properties of the extruded polymer. In general, a fairly accurate estimation of the material friction and the viscosity behavior is extremely hard to achieve in such processes due to the change in material composition, and the strong interaction between the heat and mass transfer phenomena. An interesting future work would be to consider an unknown time-dependent perturbation acting on the polymer convection speed. As it is shown in the simulation results, the state-dependent compensator might be able to handle time-varying perturbations acting on the vector field and the delay function. The proof of the robustness properties of the state-dependent input delay predictor is an interesting perspective.\\\\
\noindent
{\bf Physical definition of the parameters}\\\\
\begin{tabular}{lll}
$L=200.10^{-3}$ & $m$&\text{Extruder length}\\
$N_0=90$ & $rpm$&\text{Screw speed}\\
$B=9.3450$  $10^{-9}$ &  $m^{4}$           &\text{Geometric parameter}           \\
$F_d=-- $ &  $Kg$    $s^{-1}$        &\text{Net forward mass flow rate}   
 \\$K_d=2.45 $ $10^{-5}$ & $m^{3}$    &\text{Nozzle conductance } \\
 $\xi=10$ $10^{-3}$ & $m$     &\text{Screw Pitch }\\
 $S_{eff}=  --$ & $m^2$ &  \text{Effective area} \\
 $V_{eff}=\xi S_{eff}$  & $m^3$ &  \text{Effective volume}\\
$\eta= --$  & $Pa$    $s^{-1}$  & \text{Melt viscosity} \\ 
 $\rho_0=1240$ &  $Kg$    $m^{-3}$  &\text{PLA Melt density} \\ 
$\epsilon=--$ & $--$ &\text{Amplitude of the perturbation}\\
$\omega=-- $  & $rad.s^{-1}$  &\text{Frequency}\\
\end{tabular}

\renewcommand{\theequation}{A-\arabic{equation}}
  \setcounter{equation}{0}  
  \section*{APPENDIX}  

We recall  Lemmas 1--8  from \cite{Bekiaris13} which are applied to the nonlinear time- and state-dependent input delay system  (\ref{delay-1}) as    Lemmas \ref{lemma1}--\ref{lemma8} for the proof of  Theorem \ref{theorem-predict}.
\begin{lem}{(\textbf{Backstepping Transformation of the Actuator State)}}\label{lemma1}
The infinite dimensional backstepping transform of the actuator state given by 
\begin{eqnarray} \label{Backstepping}
W(\theta) &=&U(\theta) - \kappa(\sigma(\theta), P(\theta)),
 \end{eqnarray}
for all  $ t-D(t,x(t))\leq \theta \leq t $,  allows to transform the system  (\ref{delay-1}) with the controller  (\ref{predict-delay-5})--(\ref{control-law}) into the following target system
\begin{align}
&\dot x(t) = f\left(t,x(t), \kappa(t,x(t))+ W(t-D(t,x(t)))\right) \label {delay-target-1} \\
&W(t)=0. \label {delay-target-1bis}
 \end{align}
\end{lem}
\textbf{Proof:} The proof of Lemma \ref{lemma1} is based on a direct verification considering  $P\left(t-D(t,x(t))\right)=x(t)$ and  $\sigma \left(t-D(t,x(t))\right)=t$ in the original system (\ref{delay-1}).
\begin{lem}{\textbf{(Inverse Backstepping Transformation)}}\label{lemma-2}
 The inverse of the infinite dimensional backstepping transormation  (\ref{Backstepping}) is defined for all $ t-D(t,x(t))\leq \theta \leq t $  by
\begin{eqnarray} \label{Inv-Backstepping}
U(\theta) &=&W(\theta) + \kappa(\bar \sigma(\theta), \Pi(\theta)),
 \end{eqnarray}
with 
 \begin{align} \label{ target-predict}
\Pi(\theta) &=
 \int_{\phi(t)}^{\theta} \frac{ f\left(\bar \sigma(s),\Pi(s), \kappa (\bar \sigma(s),\Pi(s))+ W(s)\right) }{1- F(\bar \sigma(s),\Pi(s),W(s))} ds\nonumber\\
&~~+x(t)\\
 F(\bar \sigma(s)&,\Pi(s),W(s))\\
&= \f{\partial D}{\partial t} \left(\bar \sigma(s),\Pi(s) \right)+\f{\partial D}{\partial x}\left(\bar \sigma(s),\Pi(s)\right)\nonumber\\
&~~\times f\left(\bar \sigma(s),\Pi(s), \kappa(\bar \sigma(s),\Pi(s))+ W(s)\right),\\ 
\bar\sigma(\theta) &=t+ \int_{\phi(t)}^{\theta} \frac{ 1}{1-F(\bar \sigma(s),\Pi(s),W(s))} ds,\label{time-predict-target-2}\\
\phi(t)&=t-D(t,x(t)).
 \end{align}
 \end{lem}

\textbf{Proof:} Direct verification considering that $P(\theta)=\Pi(\theta)$ and $\bar \sigma(\theta)=\sigma(\theta)$ for all $ t-D(t,x(t))\leq \theta \leq t $. We refer to $P(\theta)$ as  the plant-predictor system and $\Pi(\theta)$ as  the  target-predictor system, respectively. However, they play different roles because they are driven by different inputs ($U$ versus $W$).

\begin{lem}{\textbf{(Stability of the Target System)}}\label{lemma-3}
 For any positive constant $g$, there exist a class $\mathcal{K}_{\infty}$ function $\rho_{\star}$ and a class  $\mathcal{KL}$  function $\beta$  such that for all solutions of the system satisfying  the feasibility condition  (\ref{feasibility}), the following holds:
\begin{eqnarray} \label{xi1}
\Xi(t) &\leq& \beta \left(\rho_{\star} (\Xi(t_0)),t-t_0\right),\quad t\geq t_0\\
\Xi(t)&=&  | x(t)  | +\sup_{ t-D(t,x(t))\leq \theta \leq t }| W(\theta)|\label{xi2}
 \end{eqnarray}
where
\begin{eqnarray} \label{predict-delay-6}
\rho_{\star}(s)=\frac{e^{ \frac {g}{1-c}}}{1-c}e^{(\frac {g}{1-c})( c_1+\mu_4(s))}s
 \end{eqnarray}
\end{lem}
\textbf{Proof:} Based on the input-to-state stability of  $\dot{x}=f\left(t,x,\kappa(t,x)+\omega\right)$  with respect to $\omega$, \textcolor{black}{namely, Assumption \ref{assumption-4}  },  there exist a smooth function $ S(t,x(t)):[t_0, \infty) \times \mathbb {R}^{n} \to \mathbb {R}_{+}$ and  class $ \mathcal {K}_{\infty}$ functions $\alpha_1$, $\alpha_2$,  $\alpha_3$, $\alpha_4$ such that
\begin{align}
&		\alpha_1(|x(t)|)\leq S(t,x(t))\leq \alpha_2(|x(t)|) \label{S-bound-1},\\ 
& \dot S(t,x(t))\leq -\alpha_3(|x(t)|)+\alpha_4(|W(t-D(t,x(t)))|)\label{S-deriv-bound},
 \end{align}
with
\begin{align} \label{Lyapunov-derivative-1}
  \dot S(t,x(t)) =& \frac{\partial S(t,x(t))}{\partial t} +   \frac{\partial S(t,x(t))}{\partial x}\nonumber\\ &\times f\left(t,x(t), \kappa(t,x(t))+ W(t-D(t,x(t)))\right).
 \end{align}
Let us define the  Lyapunov function  for the target system  (\ref {delay-target-1}) and  (\ref {delay-target-1bis}) as
\begin{eqnarray} \label{Lyapunov}
		V(t)=S(t,x(t))+k \int_0^{L(t)} \frac{ \alpha(r)}{r} dr,
 \end{eqnarray}
where
\begin{align} \label{L}
		L(t)&=\sup_{ t-D(t,x(t))\leq \theta \leq t }\Big | e^{g(1+\sigma(\theta)-t)}W(\theta)\Big| \nonumber \\
&= \lim_{n \to \infty}\Big(  \int_{ t-D(t,x(t))}^{t}  e^{2n g(1+\sigma(\theta)-t)} W(\theta)^{2n} d \theta  \Big)^{1/2n},
 \end{align}
with $g>0$.  Let us upperbound and lowerbound (\ref{L}) in terms of:  \be \sup_{ t-D(t,X(t))\leq \theta \leq t }\Big |W(\theta)\Big|\ee

\begin{itemize}
\item \textbf{Upperbound of $L(t)$}:
Using the feasibility condition    (\ref{feasibility}) and (\ref{sigma-deriv}), we deduce 
\begin{eqnarray} \label{predict-delay-6}
		\dot \sigma(\theta)\leq \frac{1}{1-c}.
 \end{eqnarray}
By integration of  (\ref{predict-delay-6})   on $[t-D(t,x(t))\, \theta]$ with $\sigma\left(t-D(t,x(t))\right)=t$, we derive the inequality
\begin{eqnarray} \label{Upper-bound}
1+ \sigma(\theta)-t &\leq  \frac{1}{1-c} \left(1+D(t,x(t))\right), \nonumber\\ & \forall \;\;  t-D(t,x(t))\leq \theta \leq t.
 \end{eqnarray}
From Assumption \ref{assumption-3} and (\ref{Feasibility-region-condition1}),  the following inequality holds:
\begin{align} \label{Upper-bound-2}
L(t)&\leq& e^{\frac{g}{1-c} \left(1+ c_1+\mu_4(|x|)\right)} \sup_{ t-D(t,x(t))\leq \theta \leq t } | W(\theta)|.
 \end{align}
\item \textbf{Lowerbound of $L(t)$}: Similarly, using the fact that $\sigma\left(t-D(t,x(t))\right)=t$ , with $\theta(t)$ being an increasing function, we obtain
\begin{eqnarray} \label{Lower-bound}
1+ \sigma(\theta)-t& \geq & 1,
 \end{eqnarray}
and hence
\begin{eqnarray} \label{Lower-bound-1}
L(t)&\geq& e^{g} \sup_{ t-D(t,x(t))\leq \theta \leq t } | W(\theta)| \;\;\; \nonumber\\& \forall&  t-D(t,x(t))\leq \theta \leq t.
 \end{eqnarray}
\end{itemize}
The time derivative of  (\ref{L}) is 
 \begin{align}\label{L-deriv1}	
 \dot L(t)=&\lim_{n \to \infty}\frac{1}{2n}\Big(  \int_{ \phi(t)}^{t}  e^{2n g(1+\sigma(\theta)-t)} W(\theta)^{2n} d \theta  \Big)^{1/2n-1}\nonumber\\
 &\times  \Big\{ -(1-\frac{dD(t,x(t))}{dt}) e^{2n g} W( t-D(t,x(t)))^{2n}\nonumber \\&- 2ng      \int_{ \phi(t)}^{t}  e^{2n g(1+\sigma(\theta)-t)} W(\theta)^{2n} d \theta\Big\},
  \end{align}
where  $\phi(t)=t-D(t,x(t))$.
By   (\ref{feasibility}),  it is clear that  $\frac{dD(t,x(t))}{dt}<1$ and hence
 \begin{align} \dot L(t)\leq- gL(t).\end{align}
Computing the derivative of the Lyapunov function ( \ref{Lyapunov}) as
\begin{align} \label{Lyapunov-derivative}
		\dot V(t)&= \dot S(t,x(t)) +k \dot L(t) \frac{\alpha_4(L(t))}{L(t)},
 \end{align}
we deduce 
\begin{eqnarray} \label{Lyapunov-derivative}
		\dot V(t)&\leq& \dot S(t,x(t)) -kg\alpha_4(L(t)).
 \end{eqnarray}
Using the boundness of $S(x(t))$, (\ref{S-deriv-bound}),  the following inequality holds:
\begin{align} \label{Lyapunov-derivative-2}
		\dot V(t)\leq-\alpha_3(|x(t)|)+\alpha_4(|W(t-D(t,x(t)))|)  -kg\alpha_4(L(t)).
 \end{align}
Imposing $ k=g^{-1}$, by (\ref{Lower-bound-1}), we derive the inequality
\begin{eqnarray} \label{Lyapunov-derivative-3}
		\dot V(t)&\leq& -\alpha_3(|x(t)|)  -\alpha_4(L(t)),
 \end{eqnarray}
and with  (\ref{S-bound-1}), (\ref{Lyapunov}) and  (\ref{L}), we conclude that there exists a $\mathcal{K}$  function $\gamma_1$ such that
\begin{eqnarray} \label{Lyapunov-derivative-4}
\dot V(t) \leq  - \gamma_1(V(t)). 
 \end{eqnarray}
By the comparison principle, there exists a   class $\mathcal{KL}$ function $\beta$ such that
\begin{eqnarray} \label{Lyapunov-derivative-4}
V(t) \leq   \beta(V(t_0), t-t_0).
 \end{eqnarray}
  From  (\ref{S-bound-1}) and (\ref{Lyapunov})  and the properties of class  $\mathcal{KL}$  functions, we finally get
\begin{eqnarray} \label{Lyapunov-derivative-4}
|x(t)| + L(t) \leq  \beta_1(|x(t_0)|+L(t_0) , t-t_0),
 \end{eqnarray}
$\beta_1$  being a class $\mathcal{KL}$ function. Considering the boundness of $L(t)$ defined in (\ref{Upper-bound}) and (\ref{Lower-bound-1}), the prove is achieved.
\cqfd

\begin{lem}{\textbf{(Bound of the Predictor in Terms of Actuator State)}}\label{lemma4}
There exists a class $\mathcal{K}_{\infty}$ function $\rho$ such that  for all the solutions of the system satisfying the feasibility condition (\ref{feasibility}), the following holds for all  $ t-D(t,x(t))\leq \theta \leq t $
\be\label{bound-predictor}
 |P(\theta)|\leq \rho \left(|x(t)|+  \sup_{ t-D(t,x(t))\leq s \leq t } | U(s)|\right).
\ee
\end{lem}
\textbf{Proof:} Differentiating (\ref{predict-delay-5}), we deduce the following relation  for all $t-D(t,x(t))\leq \theta \leq t$
\begin{align} \label{deriv-predict}
 \frac{dP(\theta)} {d\theta}&= \frac{ f(\sigma(\theta),P(\theta), U(\theta))}{1 -F(\sigma(\theta),P(\theta),U(\theta)) },\nonumber\\
F(\sigma(\theta),P(\theta),U(\theta))&= \f{\partial D}{\partial x} D\left(\sigma(\theta),P(\theta)\right) f(\sigma(\theta),P(\theta),U(\theta))\nonumber \\
&+ \f{\partial D}{\partial t} \left(\sigma(\theta),P(\theta) \right), 
 \end{align}
and with  the change of variable $y=\sigma(\theta)$,  (\ref{deriv-predict}) may be rewritten as:
\begin{align} \label{deriv-predict-newold0}
&\frac{dP(\phi(y))} {dy}=f\left(y,P(\phi(y)), U(y- D(\phi(y)))\right),\\ &  t\leq y \leq \sigma( t).\nonumber 
 \end{align}
From  Assumption \ref {assumption-2} we get that
\begin{align} \label{deriv-predict-newold}
\frac{dR(y,P(\phi(y)))} {d\theta}&\leq \dot \sigma(\theta) \left(R(y,P(\phi(y)))\right. \nonumber\\ 
&+ \left. \mu_3\left(| U(y- D(\phi(y)))|\right)\right),  
 \end{align}
for all $ t\leq y \leq \sigma( t)$ and using  the feasibility condition  (\ref{feasibility}), we deduce, for all $ t-D(t,x(t))\leq \theta \leq t$
\begin{align} \label{deriv-predict-new0}
& \frac{dR(\sigma(\theta),P(\theta))} {d\theta}\leq \frac{1}{1-c}\left(R(\sigma(\theta),P(\theta))+  \mu_3(|U(\theta)|)\right).
 \end{align}
By  Assumption \ref{assumption-3} and the comparison principle, we obtain
\begin{align} \label{deriv-predict-new}
 & R(\sigma(\theta),P(\theta))\leq e^{\frac{1}{1-c} \left( c_1+\mu_4(|x|)\right)} \nonumber\\ &\left(R(t,x(t)) +\sup_{ t-D(t,x(t))\leq s \leq t } \mu_3( | U(s)|)\right),\nonumber\\
 & t-D(t,x(t))\leq \theta \leq t. 
 \end{align}
With the standard properties of class $\mathcal{K}_{\infty}$ functions  the Lemma  (\ref{lemma4}) is deduced and the class $\mathcal{K}_{\infty}$ function $\rho$  is written as:
\be
\rho(s)=\mu_1^{-1}\left( \left(\mu_2(s) +\mu_3(s)\right)e^{\frac{1}{1-c} \left( c_1+\mu_4(s)\right)}  \right).
\ee
\cqfd
\begin{lem}{\textbf{(Bound of the Predictor in Terms of Transformed   Actuator State)}}\label{lemma5}
There exists a class $\mathcal{K}$ function $\psi$ such that for all the solutions of the system satisfying the feasibility condition (\ref{feasibility}), the following holds:
\be\label{bound-predictor-transform}
 |\Pi(\theta)|\leq \psi \left(|x(t)|+  \sup_{ t-D(t,x(t))\leq s \leq t } | W(s)|\right),
\ee
 \rm{for all}  $ t-D(t,x(t))\leq \theta \leq t $
\end{lem}

\textbf{Proof:} The plant $\dot X(t) = f\left(t,x(t), \kappa(t,x(t))+ \omega(t)\right) $
satisfying  the uniform  input-to-state stability property with respect to $ \omega$￼, and the function $\kappa$ ￼being locally Lipschitz in both arguments and uniformly bounded with respect to its first argument, there exist a class $\mathcal{KL}$ function $\beta_2$ and a class $\mathcal{K}$  function $\psi_1$ such that for all $\tau \geq t_0$
\begin{align} \label{target-2}
Y(\tau)&\leq& \beta_2 \left(|Y(t_0)|, \tau-t_0\right) +\psi_1 \left( \sup_{  s \geq t_0 }  | \omega(s)|\right), 
 \end{align}
with
\begin{align} \label{target-3}
\dot Y(\tau)&=& f\left(Y(\tau), \kappa(\tau,Y(\tau))+ \omega(\tau)\right). 
 \end{align}
Now, we consider the change of variable $y=\sigma(\theta)$ and write   the predictor of the target system   (\ref{ target-predict})   as
\begin{align} \label{target-4}
& \frac{d \Pi(\phi(y))} {dy}=f\left(y,\Pi(\phi(y)),  \kappa(y, \Pi(\phi(y)) )+ \omega(\phi(y))\right), \nonumber\\  & t\leq y \leq \sigma( t). 
 \end{align}
Using  (\ref{target-3}), we derive the following relation
\begin{align} \label{target-2}
\hspace{-0.5cm}&|\Pi(\theta)|\leq \psi_2 (|x(t)| )+\psi_1 \left( \sup_{ t-D(t,x(t))\leq s \leq t }  | W(s)| \right),
 \end{align}
for all  $t-D(t,x(t))\leq \theta \leq t $ with a class $\mathcal{K}$ function $ \psi_2(s)=\beta_2(s,0) $ . Using the properties of class $\mathcal{K}$ functions,  (\ref{bound-predictor-transform}) is deduced with $ \psi(s)=\psi_1 (s) +\psi_2(s)$. \cqfd 

\begin{lem}{\textbf{(Equivalence of the Norms of the Original and the target system)}}\label{lemma6}\\
There exist  class $\mathcal{K}_{\infty}$ functions $\rho_1$, $\mu_7$ such that for all the solutions of the system satisfying the feasibility condition  (\ref{feasibility}) and  for all   $t \geq t_0$, the following hold: 
\begin{eqnarray}
 \Omega(t)&\leq& \mu_7^{-1}(\Xi(t)),\label{norm-equivalence-1} \\
 \Xi(t)&\leq& \rho_1(\Omega(t)), \label{norm-equivalence-2}
\end{eqnarray}
where $\Omega$ and $\Xi$ are defined in \eqref{predict-delay-5bis-2} and \eqref{xi2}, respectively.
\end{lem}
\textbf{Proof:} 
 Using  the inverse transformation  (\ref{Inv-Backstepping})  and the bound  (\ref{bound-predictor-transform}), we derive  (\ref{norm-equivalence-1}) with
\be \mu_7^{-1}(s)=s+ \hat \rho ( \psi(s)),\ee
and from the direct  transformation  (\ref{Backstepping}) together with  the bound  (\ref{bound-predictor}), we deduce  (\ref{norm-equivalence-2}), where $\rho_1$ is define as  
\be \rho_1(s)=s+ \hat \rho( \rho(s))\ee \cqfd
\begin{lem}{\textbf{(Ball Around the Origin Within the Feasibility Region)}}\label{lemma7}
There exists a positive constant  $\bar \gamma$ such that for all the solutions of the system that satisfy
\begin{eqnarray} \label{Feasibility-region-1}
  | x(t)  | +\sup_{ t-D(t,x(t))\leq \theta \leq t }| U(\theta)| < \bar \gamma,
 \end{eqnarray} 
 the feasibility  condition  (\ref{feasibility}) is satisfied.
\end{lem}
\textbf{Proof:} \textcolor{black}{ From \eqref{f-restrict} we derive the following inequality}
\begin{eqnarray} \label{system-estimate-1}
&| f\left(t,x(t), U(t-D(t,x(t))\right)| \leq  \nonumber\\ &  \hat \alpha \left(|x(t)|  +\sup_{ t-D(t,x(t))\leq s \leq t }| U(s)| \right) 
 \end{eqnarray} 
Recalling  the relations (\ref{Feasibility-region-condition3}) and  ( \ref{Feasibility-region-condition2})  of Assumption \ref{assumption-3}, we deduce that  for all $ \theta \in [ t-D(t,x(t)), t]$ and $c \in ]0,1[$ , if a solution satisfies  
\begin{align} \label{Feasibility-region-122}
& c_3+\mu_5(|P(\theta)|)+ \left(c_2+ \mu_6(|P(\theta)|) \right )\nonumber\\ & \hat \alpha \left( | P(\theta)  | +\sup_{ t-D(t,x(t))\leq s \leq t }| U(s)| \right) < c,
 \end{align} 
 then it also satisfies (\ref{feasibility}).\\
 Using Lemma \ref{lemma4} we conclude that  \eqref{Feasibility-region-122} is satisfied if the following holds
\begin{align} \label{Feasibility-region-4}
& \left(c_2+ \mu_6(\rho(\Omega(t)))  \right ) \hat \alpha \left( \rho(\Omega(t))   +\Omega(t) \right)\nonumber\\ &+\mu_5(\rho(\Omega(t))) < c-c_3.
 \end{align} 
 Let us define  a class  $\mathcal{K}_{\infty}$ function $\rho_c$ as 
\begin{eqnarray} \label{Feasibility-region-5}
\rho_c(s) =\mu_5(\rho(s))+ \left(c_3+ \mu_6(\rho(s))  \right ) \hat \alpha \left( \rho(s)   +s \right).
 \end{eqnarray} 
It follows that
\begin{eqnarray} \label{Gamma}
\bar \gamma=\rho^{-1}_c(c-c_3).
 \end{eqnarray}  \cqfd
\begin{lem}{\textbf{(Estimate of the Region of Attraction)}}\label{lemma8}
There exists a class $\mathcal{K}$ function $\psi_{\rm RoA}$ such that for all initial conditions of the closed-loop system that satisfy relation (\ref{Feasibility-region-theo}), the solutions of the system satisfy  (\ref{Feasibility-region-1}) for $c \in ]0,1[$,  and hence, satisfy (\ref{feasibility}).
\end{lem}

\textbf{Proof:} Using Lemma \ref{lemma6} and  (\ref{xi1}), the following holds:
\begin{eqnarray} \label{Estim-1}
\Omega(t) \leq \mu_7^{-1} \left (\beta(\rho_{\star} (\rho_1(\Omega(t_0))),t-t_0)\right),
 \end{eqnarray} 
where $\Omega$ is defined in  \eqref{predict-delay-5bis-2}. Introducing the class $\mathcal{K}_{\infty}$ function  $\mu_9(s)= \mu_7^{-1} \left (\beta(s,0)\right)$, we derive the inequality
\begin{eqnarray} \label{Estim-2}
\Omega(t)  \leq \mu_9\left (\rho_{\star} (\rho_1(\Omega(t_0)))\right).
 \end{eqnarray} 
\textcolor{black}{Hence, for all initial conditions that satisfy the bound (\ref{feasibility})  with any class   $\mathcal{K}$  choice  \be \psi_{RoA}(c-c_3) \leq  \rho_1^{-1}(\rho_*^{-1}(\mu_9^{-1}(\rho_c^{-1}(c-c_3)))), \ee
the solutions satisfy (\ref{Feasibility-region-1}). 
Moreover, for all of those initial conditions, the solutions verify (\ref{predict-delay-5bis-2}),  for all $\theta > t_0-D(t_0,x(t_0))$.\cqfd}
$ $

\end{document}